\documentclass[12pt]{amsart}
\usepackage{amssymb}
\usepackage{epsf}
\usepackage{colordvi}

\numberwithin{equation}{section}
\newtheorem{thm}{Theorem}[section]
\newtheorem{defn}[thm]{Definition}
\newtheorem{prop}[thm]{Proposition}
\newtheorem{lemma}[thm]{Lemma}
\newtheorem{cor}[thm]{Corollary}
\newtheorem{example}[thm]{Example}
\newtheorem{remark}[thm]{Remark}
\newenvironment{ex}{\begin{example}\rm}{\end{example}}
\newenvironment{dfn}{\begin{defn}\rm}{\end{defn}}
\newenvironment{rem}{\begin{remark}\rm}{\end{remark}}

\newcounter{FNC}[page]
\def\fauxfootnote#1{{\addtocounter{FNC}{2}$^\fnsymbol{FNC}$%
     \let\thefootnote\relax\footnotetext{$^\fnsymbol{FNC}$#1}}}
\newcommand{\Fl}{{\mathbb{F}\ell}}

\newcommand{\Apdot}{{A_\bullet'}}
\newcommand{\Adot}{{A_\bullet}}
\newcommand{\Bdot}{{B_\bullet}}

\newcommand{\Edot}{{E_\bullet}}
\newcommand{\Epdot}{{E_\bullet'}}
\newcommand{\Fdot}{{F_\bullet}}
\newcommand{\Fpdot}{{F_\bullet'}}
\newcommand{\Gdot}{{G_\bullet}}
\newcommand{\Gpdot}{{G_\bullet'}}

\newcommand{\frakS}{{\mathfrak{S}}}
\newcommand{\calD}{{\mathcal{D}}}
\newcommand{\calE}{{\mathcal{E}}}
\newcommand{\calF}{{\mathcal{F}}}
\newcommand{\calG}{{\mathcal{G}}}
\newcommand{\calH}{{\mathcal{H}}}
\newcommand{\calI}{{\mathcal{I}}}
\newcommand{\calO}{{\mathcal{O}}}
\newcommand{\calR}{{\mathcal{R}}}
\newcommand{\calS}{{\mathcal{S}}}
\newcommand{\C}{{\mathbb{C}}}
\newcommand{\Z}{{\mathbb{Z}}}

\newcommand{\Span}[1]{{\langle #1 \rangle}}

\title[Grothendieck Polynomials]{Grothendieck Polynomials via Permutation
Patterns and Chains in the Bruhat Order}

\author{Cristian Lenart}
\address{Department of Mathematics and Statistics\\
         State University of New York at Albany\\
         Albany, NY \ 12222\\
         USA}
\email{lenart@albany.edu}
\urladdr{http://math.albany.edu/math/pers/lenart}

\author{Shawn Robinson}
\address{Department of Mathematics\\
         Rutgers University\\
         Piscataway, NJ \ \\
         USA}
\email{shawnr@math.rutgers.edu}

\author{Frank Sottile}
\address{Department of Mathematics\\
         Texas A\&M University\\
         College Station\\
         Texas \ 77843\\
         USA}
\email{sottile@math.tamu.edu}
\urladdr{http://www.math.tamu.edu/$\sim${}sottile}

\thanks{Research of Robinson supported in part by the NSF VIGRE program}
\thanks{Research of Sottile supported in part by NSF CAREER grant
  DMS-0134860, the Clay Mathematical Institute, and the MSRI}

\copyrightinfo{}{}
\begin{document}

\begin{abstract}
 We give new formulas for Grothendieck polynomials of two types.
 One type expresses any specialization of a Grothendieck polynomial in at least
 two sets of variables as a linear combination of products Grothendieck
 polynomials in each set of variables, with coefficients Schubert structure
 constants for Grothendieck polynomials.  
 The other type is in terms of chains in the Bruhat order.
 We compare this second type to other constructions of Grothendieck 
 polynomials within the more general context of double Grothendieck polynomials
 and the closely related $\calH$-polynomials. 
 Our methods are based upon the geometry of permutation patterns.
\end{abstract}

\maketitle

\tableofcontents

\section*{Introduction}

In 1982, Lascoux and Sch\"utzenberger introduced Grothendieck
polynomials~\cite{LS83}, which are inhomogeneous polynomials representing
classes of  
structure sheaves of Schubert varieties in the Grothendieck ring of the flag  
variety.
This initiated the combinatorial study of these polynomials, with Fomin 
and Kirillov giving a combinatorial construction of Grothendieck
polynomials~\cite{FK94} in terms of rc-graphs (see also~\cite{LS82b}).  
Knutson and Miller~\cite{KM03a,KM03b} gave a geometric and
commutative algebraic construction of Grothendieck polynomials, recovering the
formula of Fomin and Kirillov.  
Lascoux gave a different construction involving alternating sign
matrices~\cite{La02}. 
We give a combinatorial construction of Grothendieck polynomials
in terms of chains in the Bruhat order, which we deduce from geometric
considerations. 

This geometry is the pullback in $K$-theory of the structure sheaf of a
Schubert variety along an embedding of the form 
 \[
   \Fl \C^m \times \Fl \C^n\ \hookrightarrow\ \Fl \C^{m+n}, 
   \leqno{(1)}
  \]
where $\Fl \C^n$ is the manifold of flags in $\C^n$.  
Such maps were studied in~\cite{BS98} and  used
in a chain-theoretic construction of Schubert
polynomials~\cite{BS02}.
Billey and Braden~\cite{BB02} showed that they are the 
geometry behind pattern avoidance in singularities of 
Schubert varieties.

Our chain-theoretic construction uses a formula to express a Grothendieck
polynomial \linebreak $\calG_w(y,x_1,\dots,x_{n-1})$ with specialized variables 
as a linear combination of polynomials of the form \linebreak 
$y^j \calG_v(x_1,\dots,x_{n-1})$.  
Lascoux and Sch\"utzenberger also gave such a formula~\cite{LS82b}, 
remarking that it  ``...est l'op\'eration duale de la
multiplication par les polyn\^{o}mes {\it sp\'eciaux}, c'est-\`a-dire de 
la {\it formule de Pieri} dans l'anneau de Grothendieck, que nous 
expliciterons
ailleurs''\footnote{...this formula is the operation dual to 
multiplication by certain special classes, that is, to a Pieri-type formula in
the  Grothendieck ring of the flag variety, which we make explicit elsewhere.}.  
We make this explicit---a Pieri-type formula for multiplication 
by a particular special Grothendieck polynomial is crucial to our proof.  
We  deduce the multiplication formula from a Monk-type formula in
$K$-theory~\cite{L02}.  
We also show that the formula for  decomposing
a Grothendieck polynomial under the substitution 
$(x_1, \dots , x_n)\mapsto(x_1,\dots,x_{k-1},y,x_k,\dots,x_{n-1})$ 
can be expressed in terms of a more general Pieri-type formula in
$K$-theory~\cite{LS02}. 
   
These formulas are particular cases of a
more general formula for specializing each variable in a Grothendieck
polynomial to a variable in one of several different sets of
variables.
Such a specialized Grothendieck polynomial is
a unique $\mathbb{Z}$-linear combination of products of
Grothendieck polynomials in each set of variables.
We identify the coefficients in this linear combination as particular
Schubert structure constants by computing the pullback in the
Grothendieck ring along maps of the form~(1). 
A consequence of our general formula is the following intriguing fact:
every coefficient of a monomial in a Grothendieck polynomial is, in a
natural way, a Schubert structure constant.  
These results generalize similar results for Schubert
polynomials~\cite{BS98,BS02}. 
Different formulas for a subclass of the specializations studied here were
obtained by Buch, Kresch, Tamvakis, and Yong~\cite{BKTY03}.

The paper is organized as follows. 
Section~\ref{S:one} reviews the combinatorics of the symmetric groups,
Section~\ref{S:two} reviews the flag manifold and its Grothendieck ring, and 
Section~\ref{S:three} discusses Grothendieck polynomials.
The basis of Schubert structure sheaves is not self-dual under the intersection
pairing in the Grothendieck ring.
Dual classes are provided by structure sheaves of the boundaries of Schubert
varieties; we describe this in Section~\ref{S:three}.
These dual classes are represented by the $\calH$-polynomials of 
Lascoux and Sch\"utzenberger~\cite{LS83}.
In Section~\ref{S:four} we discuss the geometry of the pattern
map, which is the basis of our formulas.
In Section~\ref{S:five} we use a substitution formula (proved in
Section~\ref{S:seven}) to prove our main result, a
chain-theoretic construction of Grothendieck polynomials. 
In Section~\ref{S:six} we compare this construction to other 
constructions and derive similar formulas for $\calH$-polynomials,
working in the more general context of double Grothendieck and $\calH$
polynomials. 
In Section~\ref{S:seven_geom}, we compute the effect of the pattern map on the
basis of Schubert classes, and deduce a formula for the specialization of a
Grothendieck polynomial at two sets of variables.
A particular case is refined in Section~\ref{S:seven} to give a formula
for specializing a single variable in a Grothendieck polynomial.
We conclude in Section~\ref{S:eight} with a general formula to express
any specialization of a Grothendieck polynomial in any number of sets of
variables.
This last formula is used in~\cite{BSY} to show that the Buch-Fulton quiver
coefficients~\cite{BF99,Bu02} as studied in~\cite{KMS03} are 
Schubert structure constants.

\medskip
{\bf Acknowledgments.} We are grateful to Allen Knutson, 
 Ezra Miller, Igor Pak, and Jim Propp for helpful comments.

\section{The Bruhat order on the symmetric group}\label{S:one}

For an integer $n>0$, let $[n]$ denote the set $\{1,2,\dotsc,n\}$.
We write $P\sqcup Q=[n]$ to indicate that $P$ and $Q$ are disjoint and their
union equals $[n]$.
Let $\calS_n$ be the symmetric group of permutations
of the set $[n]$.
We represent permutations with one-line notation so that $w=[w_1,w_2,\dotsc,w_n]$.
Commas and brackets are sometimes omitted, so that $1432=[1,4,3,2]$.
Write $e$ for the identity element in the symmetric group, and 
$\omega_0=[n, n{-}1, \dotsc,2,1]$ for the longest element in $\calS_n$.
For statements involving more than one symmetric group, we write 
$\omega_n$ to indicate the longest element in $\calS_n$.

A permutation $w$ has a {\bf descent} at position $i$ if $w(i)>w(i+1)$.
For $a<b$ let $t_{a,b}$ be the transposition of the
numbers $a$ and $b$ and let $s_i:=t_{i,i+1}$, a simple transposition. 
The {\bf length}, $\ell(w)$, of a permutation $w\in\calS_n$ is the length of a
minimal factorization of $w$ into simple transpositions, called a 
{\bf reduced word} for $w$. 
This length is also the number of inversions of $w$ (pairs $a<b$ with $w(a)>w(b)$).
A permutation is determined by its set of inversions, and
this fact enables the following important definition.
For $w\in\calS_n$ and a set $P=\{p_1<p_2\dotsb<p_k\}\subset[n]$, let
$w|_P\in\calS_k$ be the permutation with the same inversions as the sequence 
$w_{p_1},w_{p_2},\dotsc,w_{p_k}$.
We say that $w$ {\bf contains the pattern} $w|_P$.
For $u\in\calS_m$ and $v\in\calS_n$, let $u\times v\in\calS_{m+n}$ be the
permutation that agrees with $u$ at positions in $[m]$ with
$(u\times v)|_P=v$, when $P=\{m{+}1,\dotsc,m{+}n\}$.

The {\bf Bruhat order} on $\calS_n$ is the order whose covers are 
$w\lessdot w\cdot t_{a,b}$, whenever $\ell(w)+1=\ell(w\cdot t_{a,b})$.
It has another characterization in that $u\leq v$ in the Bruhat order 
if and only if for each $i=1,\dotsc,n$, we have
$\{u_1,\dotsc,u_i\}\leq \{v_1,\dotsc,v_i\}$, where subsets of the same 
size
are compared component-wise, after putting them in order.
This {\bf tableau characterization} is implicit in the work of
Ehresmann~\cite{Eh34}. 

\section{The flag variety and its Schubert varieties}\label{S:two}

Let $V$ be an $n$-dimensional complex vector space.
Write $\Fl V$ for the set of full flags of subspaces of $V$
\[
  \Fl V\ :=\  \{ \{0\}=V_0 \subset V_1 \subset \cdots \subset V_{n-1} 
     \subset V_n = V \, | \, \dim_{\C} V_i = i, \, 0 \leq i \leq n  \}\,.
\]
The group $GL(V)$ acts transitively on $\Fl V$.
Let $E_{\bullet}$ be a fixed reference flag in $\Fl V$.  
The stabilizer of $E_{\bullet}$ is a Borel subgroup $B$, 
so that $\Fl V=GL(V)/B$.
This {\bf flag variety} $\Fl V$ is the disjoint union of $B$-orbits indexed by
elements of the symmetric group $\calS_n$ 
\[
  \Fl V\ =\ GL(V)/B\ =\ \bigsqcup_{w \in \calS_n} B \cdot wB\,.
\]
The orbit $B\cdot wB$ is isomorphic to complex affine space of dimension 
$\ell(w)$.  
The {\bf Schubert cell} indexed by $w$ is
$X_w^{\circ} = B\cdot\omega_0wB$, and the 
{\bf  Schubert variety} $X_w$ indexed by $w$ is the closure of
$X_w^{\circ}$ in  $GL(V)/B$.  
Thus $X_w$ is a subvariety of complex
codimension $\ell(w)$.  
Ehresmann~\cite{Eh34} originally defined Schubert varieties in
linear-algebraic terms:
 \begin{equation}\label{E:Ehresmann}
   X_w\ =\ X_w(E_{\bullet})\ := \
    \bigl\{ F_{\bullet} \in \Fl V \mid \dim (F_p\cap E_q) \geq 
     \#\{i \leq p \mid w(i)>n-q\} \  \forall p,q\bigr\}.
 \end{equation}

The {\bf Schubert class} $[X_w] \in H^{2\ell(w)}(\Fl V)$ is the
cohomology class Poincar\'{e} dual to the fundamental homology class of
$X_w$ in  $H_*(\Fl V),$
and these Schubert classes form a basis for the cohomology of $\Fl V$. 
If we choose a second reference flag $E'_{\bullet}$ opposite to 
$E_{\bullet}$ and set $X_v':= X_v(E'_{\bullet})$, then $[X_v] = [X_v']$ 
and the intersection of Schubert varieties corresponds to the product 
in cohomology   
\[
  [X_w \cap X'_v]\ =\ [X_w]\cdot [X_v] \in H^*(\Fl V), \  \text{ for any }
    w, \, v  \in  \calS_n. 
\]

Cohomology has an intersection pairing, defined by 
$(\alpha,\beta)\mapsto \deg(\alpha\cup\beta)$, where 
$\deg\colon H^*\to \mathbb{Z}=H^*(\textrm{pt})$ is the pushforward in cohomology
under the map to a point. 
The Schubert cohomology basis is self-dual with respect to this pairing
in that 
\[
  \deg\bigl([X_w] \cdot [X_v]\bigr)\ =\ \left\{\begin{array}
                    {c@{\quad  \quad} l}
                    1 & \text{ if } v=\omega_0w, \\
                    0         & \text { otherwise.}
                    \end{array} \right.
\]

For each $i=1,\dotsc,n$, let $V_i$ be the tautological rank $i$ 
vector bundle on $\Fl V$, $L_i = V_{n+1-i}/V_{n-i}$ be the quotient line
bundle  on $\Fl
V$, and $c_1(L_i) \in H^2 (\Fl V)$ be the first Chern class of $L_i$.  
The assignment $c_1(L_i) \mapsto x_i$ induces a ring isomorphism
\[
  H^*(\Fl V)\ \simeq\  \Lambda\ :=\ \Z[x_1, \ldots, x_n]/ I_n,
\]
where $I_n$ is the ideal generated by the nonconstant homogeneous symmetric
polynomials. 
Let $R_n(x)$ be the linear span of monomials in $\Z[x_1,\ldots, x_n]$ that
divide $x_1^{n-1}x_2^{n-2}\dotsb x_{n-1}$.
Elements of $R_n(x)$ form a complete transversal to $I_n$ and so 
the images of the monomials from $R_n(x)$ give an integral basis for $\Lambda$.

Since $\Fl V$ is smooth, the Grothendieck ring $K^0(\Fl V)$  of vector
bundles on $\Fl V$ is also the Grothendieck group of coherent 
sheaves on $\Fl V$.
A distinguished basis for $K^0(\Fl V)$ is provided by the 
{\bf Schubert classes} 
$[\calO_{X_w}]$ for $w\in\calS_n$, where $\calO_{X_w}$ is the structure 
sheaf of the Schubert variety $X_w$.
Brion showed~\cite[Lemma 2]{Br01} that the product of Schubert classes in
$K^0(\Fl V)$ corresponds to the intersection of Schubert varieties 
 \begin{equation}
   [\calO_{X_w}] \cdot [\calO_{X_v}]\ =\ 
      [\calO_{X_w \cap X'_v}] \quad \text{ for any } w, v \in \calS_n\,.
 \end{equation}

If $L_i^\vee$ is the dual of line bundle $L_i$, then the map 
$[L_i^\vee] \mapsto 1-x_i$ induces a ring isomorphism
\[
  K^0(\Fl V)\ \simeq\ \Lambda\ =\ 
    \Z[x_1, \ldots, x_n]/I_n\,,
\]
identifying the cohomology and Grothendieck rings of the flag variety.  

The intersection form in $K^0(\Fl V)$ (and hence the basis dual to the 
Schubert basis) is not as widely known as that for cohomology, so we 
review
it in more detail.
The map $\Fl V\to {\rm pt}$ induces the pushforward map on sheaves
\[
   \chi \  \colon \ \calE\ \longmapsto \ 
   \sum_{i\geq 0} (-1)^i h^i(\calE)\,,
\]
and hence a map $\chi\colon K^0(\Fl V)\to K^0({\rm pt})=\mathbb{Z}$.
Here, $h^i(\calE)$ is the dimension of the $i$th sheaf cohomology group 
$H^i(\Fl V,\calE)$ of the coherent sheaf $\calE$.
The Grothendieck ring has a natural intersection pairing given by the map
\[
   ([\calE], [\calE'])\, \longmapsto\, \chi([\calE]\cdot[\calE'])\,.
\]
Since $\chi(\calO_{X_w})=1$ for any
permutation $w$, and if $\omega_0w\geq u$, then 
$\chi(\calO_{X_w\cap X'_u})=1$,~\cite[Section 3]{BrLa01}, the 
Schubert classes do not form a self-dual basis as in cohomology.

However, if we let $\calI_w$ be the ideal sheaf of the 
boundary $X_w-X_w^{\circ}$  of the Schubert variety $X_w$, then Brion and
Lakshmibai~\cite{BrLa01} show that the classes $[\calI_w]$ form a basis dual to
the Schubert basis, 
\[
  \chi \bigl([\calO_{X_w}]\cdot[\calI_v]\bigr)\ =\ 
    \left\{\begin{array} {c@{\quad  \quad} l}
                   1 & \text{ if } v=\omega_0w, \\
                   0         & \text { otherwise.}
                   \end{array} \right.
\]
We use an expression for $[\calI_w]$ in terms of the 
Schubert classes.
We are indebted to Allen Knutson who communicated to us this expression, 
as well the proof given below.

\begin{prop}\label{P:duality}
 ${\displaystyle [\calI_w]\ =\ \sum_{v\geq 
w}(-1)^{\ell(vw)}[\calO_{X_v}]}$.
\end{prop}

\begin{proof}
 We show that the expression on the right hand side defines elements of 
 $K^0(\Fl V)$ that are dual to the Schubert classes.
 Consider
 \begin{eqnarray*}
  \chi \Bigl( [\calO_{X_u}]\cdot \sum_{v\geq w}
                 (-1)^{\ell(vw)}[\calO_{X_v}]\Bigr)
    &=& \sum_{v\geq w}(-1)^{\ell(vw)}
                  \chi\bigl( [\calO_{X_u}]\cdot[\calO_{X_v}]\bigr) \\
    &=& \sum_{v\geq w}(-1)^{\ell(vw)}
                  \chi\bigl([\calO_{X_u\cap X'_v}]\bigr) \\
    &=& \sum_{\omega_0u\geq v\geq w}(-1)^{\ell(vw)}
                  \chi\bigl( [\calO_{X_u\cap X'_v}]\bigr)\,.
 \end{eqnarray*}
 The restriction in the index of summation of this last sum is because 
 $X_u\cap X'_v=\emptyset$ unless $\omega_0u\geq v$.
 Since
 $\chi\bigl( [\calO_{X_u\cap X'_v}]\bigr)=1$ for such $v$, this last sum is
\[
  \sum_{\omega_0u\geq v\geq w}(-1)^{\ell(vw)}\ =\ 
\left\{\begin{array}{rcl}
     1&\ &\mbox{if }\omega_0u=w\\
     0&&\mbox{otherwise}\end{array}\right.
\]
 as $(-1)^{\ell(vw)}$ is the M\"obius function of the
 interval $[w,v]$ in the Bruhat order~\cite{De77}. 
\end{proof}

Since the Schubert classes  $[\calO_{X_w}]$ form a basis for 
$K^0(\Fl V)$, there are integral {\bf Schubert structure constants}
$c^w_{u,v}$ for $u,v,w\in\calS_n$ defined by the identity
\[
  [\calO_{X_u}]\cdot[\calO_{X_v}]\ =\ 
  \sum_w c^w_{u,v}\,[\calO_{X_w}]\,.
\]
Since $[\calI_{\omega_0w}]$ is dual to $[\calO_{X_w}]$, we also have
\begin{equation}\label{E:StrConst}
   c^w_{u,v}\ =\
   \chi ([\calO_{X_u}]\cdot[\calO_{X_v}]\cdot[\calI_{\omega_0w}])\,.
\end{equation}
The only general formulas known for the Schubert structure constants are for the  
multiplication by $[\calO_{X_{s_i}}]$, the class of a hypersurface Schubert
variety or of an ample line bundle~\cite{Br01,FL94,La90b,L02,LP03,LiSe03,PR99},
and a Pieri-type formula for multiplication by certain
special classes in $K^0(\Fl V)$~\cite{LS02}. 
In Section~\ref{S:five} we use Lenart's formula~\cite{L02} for multiplication of
a Schubert class in $K^0(\Fl V)$ by the class $[\calO_{X_{s_1}}]$, which is
represented by the variable $x_1$. 

\section{Schubert and Grothendieck polynomials}\label{S:three}

Lascoux and Sch\"utzenberger~\cite{LS82a} introduced 
{\bf Schubert polynomials}, which are polynomials in $\mathbb{Z}[x_1,\dotsc,x_n]$
representing the Schubert cohomology classes, and later~\cite{LS83} 
{\bf Grothendieck polynomials},  which are polynomial representatives of 
the Schubert classes in the Grothendieck ring.
Elements of the symmetric group $\calS_n$ act on polynomials in
$\mathbb{Z}[x_1,\dotsc,x_n]$ by permuting the indices of variables.

For each $i=1,\dotsc,n{-}1$, define operators $\partial_i$ and $\pi_i$
on $\mathbb{Z}[x_1,\dotsc,x_n]$:  
\begin{equation}\label{def:pi}
   \partial_i := \frac{1-s_i}{x_i -x_{i+1}} \qquad \pi_i := 
   \partial_i(1-x_{i+1})\,.
\end{equation}
Suppose $w \in \calS_n$ has reduced decomposition 
$w=s_{i_1} \cdots s_{i_r}$.  
The operators $\partial_w := \partial_{i_ 1}\cdots\partial_{i_r}$ and 
$\pi_w := \pi_{i_1} \cdots \pi_{i_r}$ are 
independent of the choice of reduced decomposition.

Set $\frakS_{\omega_0} = \calG_{\omega_0} = x_1^{n-1}x_2^{n-2}\cdots 
x_{n-1}$.  
The Schubert polynomial representing $[X_w]$  is 
$\frakS_w = \partial_{w^{-1}\omega_0}\frakS_{\omega_0}$, and the 
Grothendieck polynomial representing $[\calO_{X_w}]$
is $\calG_w = \pi_{w^{-1}\omega_0} \calG_{\omega_0}$.  
The Schubert polynomial $\frakS_w$ is the lowest degree
homogeneous part of the Grothendieck polynomial $\calG_w$. 

\begin{ex}\label{E:S3}
  For example, we give the Grothendieck polynomials $\calG_w$ for $w\in\calS_3$:
\[
  \begin{picture}(263,90)(0,-3)
   \put( 60, 75){$\calG_{321}(x)\ =\ x_1^2x_2$}
   \put(130, 50){$\calG_{231}(x)\ =\ x_1x_2$}
   \put(  0, 50){$\calG_{312}(x)\ =\ x_1^2$}
   \put(  0, 25){$\calG_{213}(x)\ =\ x_1$}
   \put(130, 25){$\calG_{132}(x)\ =\ x_1+x_2-x_1x_2$}
   \put( 60,  0){$\calG_{123}(x)\ =\ 1$}
  \end{picture}
\]
 Also, in the special case when when the permutation $w$ has a unique
 descent at $1$ so that $w=j{+}1,1,\dotsc,j,j{+}2,\dotsc,n$, then 
\[
  \calG_{j{+}1,1,\dotsc,j,j{+}2,\dotsc,n}(x)\ =\ x_1^j\,.
\]
\end{ex}

The sets $\{\frakS_w\mid w\in\calS_n\}$ and  
$\{\calG_w\mid w\in\calS_n\}$ each form bases for  
$R_n(x)$. 
In particular, if $f\in R_n(x)$ and and $g\in R_n(x)$, then
$f\cdot g\in R_{m+n}(x)$.
If we specialize the variables in $f\in R_n(x)$ to two
different sets of variables $y$ and $z$ as in 
the map $\psi_{P,Q}$ of Section~\ref{S:four}, then the
resulting polynomial lies in $R_n(y)\otimes R_n(z)$.
The identifications of $R_n(x)$ with the cohomology and Grothendieck 
rings enable us to use formulas in the cohomology and Grothendieck rings
to deduce formulas for the polynomials $\frakS_w$ and $\calG_w$.

The symmetric group $\calS_m$ naturally embeds into any symmetric group
$\calS_n$ for $m<n$, and so we set $\calS_\infty:=\cup\calS_n$.
The polynomials $\frakS_w$ and $\calG_w$ do not depend upon $n$ as long
as $w\in\calS_n$, and thus are well-defined for $w\in\calS_\infty$.
This fact allows us to use geometry (which {\it a priori} concerns Schubert
classes in $K$-theory) to establish identities of polynomials.
The sets
\[
   \{\frakS_w\mid w\in\calS_\infty\}\qquad\mbox{and}\qquad
   \{\calG_w\mid w\in\calS_\infty\}
\]
form integral bases for the polynomial ring
$\mathbb{Z}[x_1,x_2,\dotsc]$.

Lascoux and Sch\"utzenberger~\cite{LS83} introduced the $\calH$-polynomials
defined for $w\in\calS_n$ by 
\begin{equation}\label{def:hpoly}
   \calH_w\ :=\ \sum_{\omega_n\geq v\geq w} (-1)^{\ell(vw)}\calG_v\,.
\end{equation}
By Proposition~\ref{P:duality}, $\calH_w$ represents the class $[\calI_w]$ in 
$K^0(\Fl\C^n)$.
The $\calH$-polynomials depend upon $n$ and
are not stable under the inclusion $\calS_n\hookrightarrow\calS_{n+1}$.

Like the Grothendieck polynomials, they have a definition in terms of
operators~\cite{La90b}.
Set $\mu_i:=\pi_i-1$, and 
$\mu_w := \mu_{i_1} \cdots \mu_{i_r}$ when 
$w=s_{i_1} \cdots s_{i_r}$ is a reduced decomposition of the permutation $w$.
If we set $\calH_{\omega_0}:=x_1^{n-1}\cdot x_2^{n-2}\dotsb x_{n-1}$, then 
$\calH_w:=\mu_{w^{-1} \omega_0}(\calH_{\omega_0})$. 
We refer to Section \ref{S:six} for various constructions of the
$\calH$-polynomials, including more details on the ones above.

\section{Geometry of the pattern map}\label{S:four}

We review the geometry of the pattern map and deduce its effects on
the monomial basis $R_n(x)$ of the Grothendieck ring.
This pattern map was studied both by Bergeron and Sottile~\cite{BS98} and 
by Billey and Braden~\cite{BB02}, but from different viewpoints. 

Let $T_0\simeq \mathbb{C}^\times$ be a 1-dimensional torus acting on a 
vector space with only two weight spaces $U$ and
$V$ of respective dimensions $m$ and $n$.
Then the vector space has an equivariant decomposition $U\oplus V$. 
Set $G:=\textit{GL}(U\oplus V)$ and let $G'$ be the centralizer of the 
torus $T_0$---a connected, reductive subgroup of $G$.
In fact, $G'$ is the group of linear transformations that respect the
decomposition $U\oplus V$, that is 
$G'=\textit{GL}(U)\times\textit{GL}(V)$.

Let $\mathcal{F}$ be the flag variety of $G$. 
The fixed points $\mathcal{F}^{T_0}$ of $\mathcal{F}$ under $T_0$ consist 
of those Borel subgroups of $G$ that contain $T_0$.
The intersection of $G'$ with such a Borel subgroup $B$ of $G$ is a Borel
subgroup $\pi(B)$ of $G'$.
This defines a map
\[
  \pi\ \colon\ \mathcal{F}^{T_0}\ \longrightarrow\ \mathcal{F}'\,.
\]
Billey and Braden~\cite[Theorem~10]{BB02} show that this map is an 
isomorphism
between each connected component of $\mathcal{F}^{T_0}$ and the flag
variety $\mathcal{F}'$ of $G'$.

Let $T\subset G$ be a maximal torus containing $T_0$.
The Weyl group $W$ of $G$ is identified with the $T$-fixed points
$\mathcal{F}^T$, which are contained in $\mathcal{F}^{T_0}$.
The torus $T$ is also a maximal torus of $G'$ and the Weyl group $W'$ of 
$G'$ is similarly identified with the $T$-fixed points $\mathcal{F}'\,^T$ of 
its flag variety.
Furthermore, we have that $\pi(W)=W'$.
 
Billey and Braden identify this map 
$\pi\colon W \to W'$ as the pattern map of Billey and 
Postnikov~\cite{BiPo}.
Briefly, they fix the maximal torus $T$ and a Borel subgroup
containing $T$, and then select $T_0\subset T$.
Then the Weyl group $W'$ is conjugate to a parabolic subgroup
$\calS_m\times\calS_n$ of $W=\calS_{m+n}$.
We describe this more explicitly below.

We interpret this map differently.
First, this same geometry was studied by Bergeron and Sottile~\cite{BS98}, 
who in addition identified the components of $\mathcal{F}^{T_0}$, as well 
as the inverse images of Schubert varieties of $\mathcal{F}'$ in each component
of $\calF^{T_0}$ as  Richardson varieties in $\mathcal{F}$.
Identifying $\mathcal{F}' = \Fl U \times\Fl V$,
then the map $\pi\colon\mathcal{F}^{T_0}\to\Fl U \times\Fl V$
is described as
 \begin{equation}\label{E:restriction}
   \Gdot \in \mathcal{F}\ \longmapsto\ 
    (\Gdot\cap U, \; \Gdot\cap V)\,.
 \end{equation}

This map is not injective.
For example, suppose that $U\simeq V\simeq\mathbb{C}^2$ with 
basis $\{\Blue{u_1},\Blue{u_2}\}$ for $U$ and 
$\{\Red{v_1},\Red{v_2}\}$ for $V$.
Then each of the flags in $\mathcal{F}^{T_0}$
 \begin{eqnarray*}
  &\Span{\Blue{u_1}}\subset\Span{\Blue{u_1},\Blue{u_2}}\subset
        \Span{\Blue{u_1},\Blue{u_2},\Red{v_1}}
        \subset\Span{\Blue{u_1},\Blue{u_2},\Red{v_1},\Red{v_2}}&\\
  &\Span{\Blue{u_1}}\subset\Span{\Blue{u_1},\Red{v_1}}\subset
        \Span{\Blue{u_1},\Red{v_1},\Blue{u_2}}
        \subset\Span{\Blue{u_1},\Red{v_1},\Blue{u_2},\Red{v_2}}&\\
  &\Span{\Red{v_1}}\subset\Span{\Red{v_1},\Blue{u_1}}\subset
        \Span{\Red{v_1},\Blue{u_1},\Blue{u_2}}
        \subset\Span{\Red{v_1},\Blue{u_1},\Red{v_2},\Blue{u_2}}&
 \end{eqnarray*}
as well as three others in $\mathcal{F}^{T_0}$ restrict to the same pair 
of flags 
\[ 
   \bigl(\Span{\Blue{u_1}}\subset\Span{\Blue{u_1},\Blue{u_2}},\ 
     \Span{\Red{v_1}}\subset\Span{\Red{v_1},\Red{v_2}}\bigr)
\]
in $\Fl U\times\Fl V$.
These six flags in $\mathcal{F}^{T_0}$ correspond to the cosets of 
$\calS_2\times\calS_2\simeq W'$ in $\calS_4=W$.

Associated to a flag $\Gdot\in\mathcal{F}$ are two subsets
 \begin{eqnarray*}
  P(\Gdot)&:=&\{j\mid \dim G_{j-1}\cap U < \dim G_j\cap U \}\ \in\ 
             {\textstyle\binom{[m+n]}{m}}\\
  Q(\Gdot)&:=&\{j\mid \dim G_{j-1}\cap V < \dim G_j\cap V \}\ \in\ 
             {\textstyle\binom{[m+n]}{n}}
 \end{eqnarray*}
Furthermore, a flag $\Gdot$ is in $\mathcal{F}^{T_0}$ if and only if these 
two subsets are disjoint.
Indeed, the flag $\Gdot$ is fixed by $T_0$ if and only if each subspace 
$G_j$ is a direct sum of $T_0$-fixed subspaces, that is 
$G_j=(G_j\cap U)\oplus (G_j\cap V)$, and this implies that $P$ and $Q$ are
disjoint.

Suppose that we have disjoint subsets $P\in\binom{[m+n]}{m}$ 
and  $Q\in\binom{[m+n]}{n}$, written as
 \begin{equation}\label{E:PQ}
   P\ =\ \{p_1<p_2<\dotsb<p_m\}\qquad\textrm{and}\qquad
   Q\ =\ \{q_1<q_2<\dotsb<q_n\}\,.
 \end{equation}
Bergeron and Sottile~\cite[\S 4.5]{BS98} study the map
\[
  \varphi_{P,Q}\ \colon\ \Fl U \times\Fl V\ \longrightarrow\ 
   \mathcal{F}^{T_0}
\]
defined by setting $\varphi_{P,Q}(\Edot,\Fdot)$ to be the flag $\Gdot$ 
such that 
$\Gdot\cap U=\Edot$ and $\Gdot\cap V=\Fdot$ with 
$P=P(\Gdot)$ and $Q=Q(\Gdot)$.
This is a section of the pattern map $\pi$.

\begin{rem}
 Lemma 4.5.1 of~\cite{BS98} shows that 
 if  $w\in\calS_n$ and $P\sqcup Q=[n]$, then there is a section $\varphi$ of the  
 pattern map such that 
\[ 
  \varphi\bigl(X_u\times X_v\bigr)\ \subset\  X_w\,,
\]
 where $u:=w|_P$ and $v:=w|_Q$, and the Schubert varieties are defined with
 respect to appropriate flags.
 (This also follows from the work of Billey and Braden.)

 This section $\varphi$ realizes
 $\varphi\bigl( X_u\times X_v\bigr)$ as a component of the $T_0$-fixed
 points of the Schubert variety $X_w$ that is isomorphic to $X_u\times X_v$.
 Thus if $X_w$ is smooth, then so are $X_u$ and $X_v$.
 This observation is due to Victor Guillemin, and is stated in the following
 proposition.
 This provides an elementary explanation for the occurrence of permutation
 patterns in the study of singularities of Schubert varieties.
 We are indebted to Ezra Miller, who communicated this to us.

\begin{prop}
 If a Schubert variety $X_w$ is smooth with $w\in\calS_n$, then so is $X_u$, 
 where $u=w|_P$, for any $P\subset[n]$.
 Conversely, if $X_u$ is singular, then so is $X_w$, for any $w$ containing $u$
 as a pattern.
\end{prop}

 Since these results (the section of the pattern map) hold for any flag variety
 $G/B$, these ideas lead to the same elementary characterization of smoothness
 via patterns. 
\end{rem}

We determine the effect of the map $\varphi_{P,Q}^*$ 
on the Grothendieck rings of these flag varieties in terms of their
monomial bases.
Suppose that $L_1,\dotsc,L_{m+n}$ are the tautological line bundles on
$\mathcal{F}$ and $M_1,\dotsc,M_m$ and $N_1,\dotsc,N_n$ are the 
tautological
line bundles in $\Fl U$ and $\Fl V$, respectively.
We pullback the bundles $M_i$ and $N_j$ to the product 
$\Fl U \times\Fl V$.
Writing $P$ and $Q$ as in~\eqref{E:PQ}, we see that 
\[
   \varphi_{P,Q}^*(L_j)\ =\ 
    \left\{\begin{array}
                   {r@{\quad} l}
                   M_i& \text{ if }\ j=p_i\,, \\
                   N_i& \text{ if }\ j=q_i\,.
                   \end{array} \right.
\]
That is, 
\[
   \varphi_{P,Q}^*(L_{p_i})\ =\ M_i\qquad\textrm{and}\qquad
   \varphi_{P,Q}^*(L_{q_i})\ =\ N_i\,.
\]
Consider the isomorphism
\[
  \psi_{P,Q}\ \colon\ \mathbb{Z}[x_1,\dotsc,x_{m+n}]\ \longrightarrow\ 
      \mathbb{Z}[y_1,\dotsc,y_m]\, \otimes\,
      \mathbb{Z}[z_1,\dotsc,z_n]
\]
 defined by
\[
   x_{p_i}\ \longmapsto\  y_i 
     \qquad\textrm{and}\qquad
   x_{q_j}\ \longmapsto\  z_j\ . 
\]
This induces a surjective algebra map
\[
  \overline{\psi}_{P,Q}\ \colon\ 
\mathbb{Z}[x_1,\dotsc,x_{m+n}]/I_{m+n}\ 
      \longrightarrow\ 
      \mathbb{Z}[y_1,\dotsc,y_m]/I_m
            \otimes
      \mathbb{Z}[z_1,\dotsc,z_n]/I_n\,.
\]

\begin{prop}\label{P:pattern_map}
  Let $P$ and $Q$ be as above.
  Then, $\varphi_{P,Q}^*=\overline{\psi}_{P,Q}$, under the identifications 
  of cohomology and Grothendieck rings with these quotients.  
\end{prop}

We will determine the effect of the map $\varphi_{P,Q}^*$ on the Schubert
basis in Section~\ref{S:seven_geom}.
In the next section we use those results to give a new construction of
Grothendieck polynomials. 

\section{A chain-theoretic construction of Grothendieck 
polynomials}\label{S:five}

We will apply Theorem~\ref{T:SingleSubstitution} of Section~\ref{S:seven}
to give an interpretation of the monomials that appear in a
Grothendieck polynomial.  Toward this end we introduce some terminology
and notation regarding the Bruhat order on the symmetric group. 

 For any $1\leq k< n$, the {\bf $k$-Bruhat order} on $\calS_n$, denoted 
 by $<_k$, is defined by its covers. 
 We say that $w \lessdot_k w\cdot t_{a,b}$ is a cover in the $k$-Bruhat order 
 if $a \leq k<b$ and if $w\lessdot w\cdot t_{a,b}$ is a cover in the
 usual Bruhat order (so that $\ell(w)+1=\ell(w\cdot t_{a,b})$).
 We label the edges of the Hasse diagram of
 $\calS_n$ by writing  
\[
   w \xrightarrow{\ (a,b)\ } w_1\quad \text{ if  }\quad w \lessdot
   w_1=w\cdot t_{a,b}.
\]
 A {\bf chain} $\gamma$ in the Bruhat order is a sequence of covers
 \[
   w \xrightarrow{\ (a_1,b_1)\ } w_1 
     \xrightarrow{\ (a_2,b_2)\ } w_2  \xrightarrow{\ (a_3,b_3)\ } \ 
\dotsb\  
     \xrightarrow{\ (a_j,b_j)\ } w_j,
\]
and we write $w(\gamma)$ for the terminal permutation in the chain.  
The  {\bf length} of the chain is
$\ell(\gamma)=\ell(w(\gamma))-\ell(w)$, the number of  transpositions in 
the chain.
A chain
\[
   w \xrightarrow{\ (a_1,b_1)\ } w_1 
     \xrightarrow{\ (a_2,b_2)\ } w_2  \xrightarrow{\ (a_3,b_3)\ } \ 
\dotsb\  
     \xrightarrow{\ (a_j,b_j)\ } w_j
\]
in the $k$-Bruhat order is {\bf increasing} if $a_i < a_{i+1}$, or 
$a_i=a_{i+1}$ and 
$b_i>b_{i+1}$ for all $1\leq i <j$.  
Any chain in the $k$-Bruhat order may be factored
uniquely as a concatenation of increasing chains, with one descent
between adjacent increasing segments.  

A {\bf marked chain} $\gamma$ in the $k$-Bruhat order is a chain in the
$k$-Bruhat order together with certain marked edges that necessarily include
the first edge. 
We require that cutting $\gamma$ at the permutation before each marked edge
decomposes it into increasing segments.
Write $m(\gamma)$ for the the number of marks on a chain $\gamma$ .
Here is a marked chain in the 1-Bruhat order
\[
  132654\ \xrightarrow{\ \underline{(1,3)}\ }\ 
  231654\ \xrightarrow{\ (1,2)\ }\ 
  321654\ \xrightarrow{\ \underline{(1,6)}\ }\ 
  421653\ \xrightarrow{\ (1,5)\ }\ 
  521643\ \xrightarrow{\ (1,4)\ }\   621543\,.
\]

We use the following special case of the Monk formula in $K$-theory.

\begin{prop}[Theorem 3.1 of~\cite{L02}] \label{T:monk-fomula}
  Let $w\in \calS_n$.  
  Then
 \begin{equation*}
   x_1 \calG_w\ =\ \sum (-1)^{l(\gamma){-}1} \calG_{w(\gamma)}\,,
 \end{equation*}
 the sum over all increasing chains $\gamma$ in the $1$-Bruhat order 
 on $\calS_n$ that begin at $w$.
\end{prop}

The following corollary is crucial to what follows.

\begin{cor}\label{C:iterated-monk}
 For any $j$, $1\leq j \leq n$, and any $w\in \calS_n$,
 \begin{equation*}
   x_1^j \calG_w = \sum (-1)^{l(\gamma){-}j} \calG_{w(\gamma)},
 \end{equation*}
 the sum  over all marked chains $\gamma$ in the $1$-Bruhat order 
 on $\calS_n$ that begin at $w$ and have $j$ marks.
\end{cor}

Corollary~\ref{C:iterated-monk} is a special case of a Pieri-type formula 
for the Grothendieck ring~\cite{LS02}.  
That formula gives the Schubert structure constants for 
multiplication by a {\bf special Schubert class},
which is indexed by a $p$-cycle
 \[
   r_{[k,p]}\ :=\
   [1,\ldots,k{-}1,k{+}p,k,k{+1},\ldots,\widehat{k{+}p},\ldots,m]\in\calS_m\,,
 \]  
where $k$ and $p$ are positive integers such that $p\leq m-k$.  
In terms of the simple generators of $\calS_m$, 
$r_{[k,p]}=s_{k+p-1} s_{k+p-2} \cdots s_k$.
For Corollary~\ref{C:iterated-monk}, $\calG_{r_{[1,j]}}=x_1^j$. 

The main result of this section requires the following lemma.
For a permutation $v\in\calS_{n-1}$, let
$n.v\in\calS_n$ be the permutation $[n, v_1, v_2, \dotsc, v_{n-1}]$.

\begin{lemma}\label{L:LR-identity}
Let $v \in \calS_{n-1}$, $w \in \calS_n$, and $0<j<n$.  Then we have the
following equality of structure constants:
 \begin{equation*}
  c_{w,r_{[1,n-1]}}^{(n-1+j).v} = c_{w,r_{[1,n-j]}}^{n.v}\,,
 \end{equation*}
where all permutations are considered as elements of $\calS_{2n-1}$.
\end{lemma}
\begin{proof}
Up to a sign, $c_{w,r_{[1,n-1]}}^{(n-1+j).v}$ is 
the number of marked chains in the 1-Bruhat order from $w$ to $(n{-}1{+}j).v$
with $n{-}1$ marks, and  $c_{w,r_{[1,n-j]}}^{n.v}$ is the number of marked
chains in the  1-Bruhat order from $w$ to $n.v$ with $n-j$ marks, 
by Corollary~\ref{C:iterated-monk}.
Since 
\[
   \ell((n-1+j).v)\ -\ \ell(n.v)\ =\ 
  \ell(r_{[1,n-1]}) \ -\ \ell(r_{[1,n-j]})\,,
\] 
the signs of these Schubert structure
constants are equal.
We will prove the  lemma by establishing a bijection between 
these two sets of chains.  

Let $\gamma$  be a chain in the 1-Bruhat order on $\calS_{2n-1}$
starting at $w$ of length at least $n{-}1$.
As $w\in \calS_n$, this chain necessarily contains an element of the form
$n.v$, for $v\in \calS_{n-1}$, and subsequent links will be 
\[
   n.v\ \xrightarrow{\ (1,n+1)\ }\ (n+1).v\ \xrightarrow{\ (1,n+2)\ }\ 
   (n+2).v\ \xrightarrow{\ (1,n+3)\ }\ \dotsb\ 
   \xrightarrow{\ (1,m)\ }\ m.v\,.
\]
Any valid marking of $\gamma$ must include these subsequent links.

The bijection is now clear.  
Given such a marked chain $\gamma$ with $m=n+1-j$ and $n{-}1$ marks, removing 
these last $j{-}1$ marked links gives a chain $\gamma'$ from $w$ to 
$n.v$ with $n{-}j$ marks.
Adding this necessary chain of length $j{-}1$ to a chain in the 1-Bruhat order
from $w$ to  $n.v$ with $n{-}j$ marks gives a chain from $w$ to 
$(n+1-j).v$ with $n{-}1$ marks.
\end{proof}

The results of Section 8 give the 
following formula which is the key step in describing the monomials that appear  
in a Grothendieck polynomial.

\begin{thm}\label{T:substitution-formula}
 Let $w\in\calS_n$.
 If $w$ has the form $n.v$, then 
 $\calG_w(y,\, x_1,\dotsc,x_{n-1})=y^{n-1}\cdot \calG_v$. 
 Otherwise, 
 \begin{eqnarray*}
   \calG_w(y,\, x_1,\dotsc,x_{n-1})
  &=& \sum_{j<n, \,  v\in S_{n-1}} c_{w,r_{[1,n-j]}}^{n.v}\, y^{j-1} \cdot 
\calG_v(x),\\
  &=& \sum_{\gamma} (-1)^{\ell(\gamma)-m(\gamma)} y^{n{-}1-m(\gamma)}\cdot 
\calG_{v}(x)\,,
 \end{eqnarray*}
 where the second sum is over all marked chains $\gamma$ in the
 $1$-Bruhat order on  $\calS_n$ that begin at $w$ and end in a
 permutation of the form $n.v$, where $v\in\calS_{n-1}$.  
\end{thm}
\begin{proof}

By Theorem~\ref{T:SingleSubstitution}, we have
\[
  \calG_w(y, x_1,\dotsc,x_{n-1})\ =\ 
   \sum_{\substack{j\geqslant 0\\v\in\calS_{n-1}}} 
     c^{n{+}j.v}_{w,r_{[1,n-1]}} y^j 
\calG_v(x)\,,
\]
and we apply Lemma~\ref{L:LR-identity} to obtain the first line.
(We caution that there is a shift of the index $j$.)
For the second, we simply use the description of 
$c_{w,r_{[1,n-j]}}^{n.v}$ given by Corollary~\ref{C:iterated-monk}.
\end{proof}

\begin{ex}\label{Ex:otherSubst}

 We may illustrate Theorem~\ref{T:substitution-formula} on $\calG_{1432}$ by
 considering only the initial segments of the chains above 1432.
 The Grothendieck polynomial $\calG_{1432}$ equals
 \begin{equation}\label{E:G1432}
  x_1^2x_2+ x_1^2x_3 + x_1x_2^2 +  x_1x_2x_3 +  x_2^2x_3
  \,- 2x_1^2x_2x_3 - x_1^2x_2^2 - 2x_1x_2^2x_3 \,+ x_1^2x_2^2x_3\,.
 \end{equation}
 Then, $\calG_{1432}(y, x_1,x_2,x_3)$ is
\[
   y^2\cdot(x_1+x_2-x_1x_2) + (y-y^2)\cdot x_1^2 + 
   (y-y^2)\cdot x_1x_2 + (1-2y+y^2)\cdot x_1^2x_2\,.
\]
 Using the list of Grothendieck polynomials in Example~\ref{E:S3}, we see that 
 this equals
\[
   y^2\cdot \calG_{132}(x) + (y-y^2)\cdot \calG_{312}(x) + 
   (y-y^2)\cdot \calG_{231}(x) + (1-2y+y^2)\cdot\calG_{321}(x)\,.
\]
 
 By Theorem~\ref{T:substitution-formula}, each Grothendieck polynomial
 $\calG_v(x)$ in this sum corresponds to a chain in the 1-Bruhat order from
 1432 to $4.v$ and the coefficient (a polynomial in $y$) corresponds
 to the different ways that such a chain may be marked.
 We write the relevant chains with the corresponding terms and leave
 possible markings to the reader.
\[
 \begin{array}{rll}
 &  1432\Blue{\xrightarrow{\,(1,2)\,}}4132
   & y^2\cdot \calG_{132}(x)\\
 &  1432\Blue{\xrightarrow{\,(1,3)\,}}
   3412\Blue{\xrightarrow{\,(1,2)\,}}4312
   & (y-y^2)\cdot \calG_{312}(x)\\
 &  1432\Blue{\xrightarrow{\,(1,4)\,}}
   2431\Blue{\xrightarrow{\,(1,2)\,}}4231
   & (y-y^2)\cdot \calG_{231}(x)\\
 &  1432\Blue{\xrightarrow{\,(1,4)\,}}
   2431\Blue{\xrightarrow{\,(1,3)\,}}
   3421\Blue{\xrightarrow{\,(1,2)\,}}4321\ \ 
   & (1-2y+y^2)\cdot\calG_{321}(x)
 \end{array}
\]

\end{ex}

Applying Theorem~\ref{T:substitution-formula} iteratively 
decomposes a Grothendieck polynomial into a sum of monomials; these are 
indexed by concatenations, left to right, of marked chains 
in the $1$--,$2$--,\ldots, $(n{-}1)$--Bruhat orders, where the permutation at
the end  of the chain in the $k$-Bruhat order has the form
$[n,n-1,\ldots,n-k+1,\ldots]$.   
We call such a chain a {\bf climbing marked chain}. 
The {\bf weight} 
$\alpha(\gamma)$ of a climbing chain $\gamma$ from $w$  to $\omega_0$ in 
the Bruhat order on $\calS_n$  is the sequence of integers 
$\alpha=\alpha_1,\alpha_2,\dotsc,\alpha_{n-1}$, where $\alpha_k$ is the 
number of markings of $\gamma$ in the segment
corresponding to the $k$-Bruhat order.  
(This segment is unique, by the condition on the form of its endpoint.)
We deduce the main result of this section, which is proved by simply iterating
Theorem~\ref{T:substitution-formula}.

\begin{thm}\label{T:monomial-formula}
 Let $w\in\calS_n$, and let $\delta = (n-1, n-2, \dots, 1).$
 Then
\[
  \calG_w\ =\ \sum (-1)^{\ell(\gamma)-m(\gamma)}
                      \, \, x^\delta/x^{\alpha(\gamma)}\,,
\]
 the sum over all climbing marked chains $\gamma$ in the Bruhat order that climb  
 from $w$ to $\omega_0$.
\end{thm}

\begin{ex}\label{fivechains}
 Let us examine the statement of 
 Theorem~\ref{T:monomial-formula} in
 a particular example.
 The Grothendieck polynomial $\calG_{1432}$ is given in~\eqref{E:G1432}. 
 The diagram below gives all increasing chains in the Bruhat order on
 $\calS_4$ from $1432$ to $4321=\omega_0$, where we have labeled each 
 cover with the corresponding transposition:
\[
   \epsfxsize=3in \epsffile{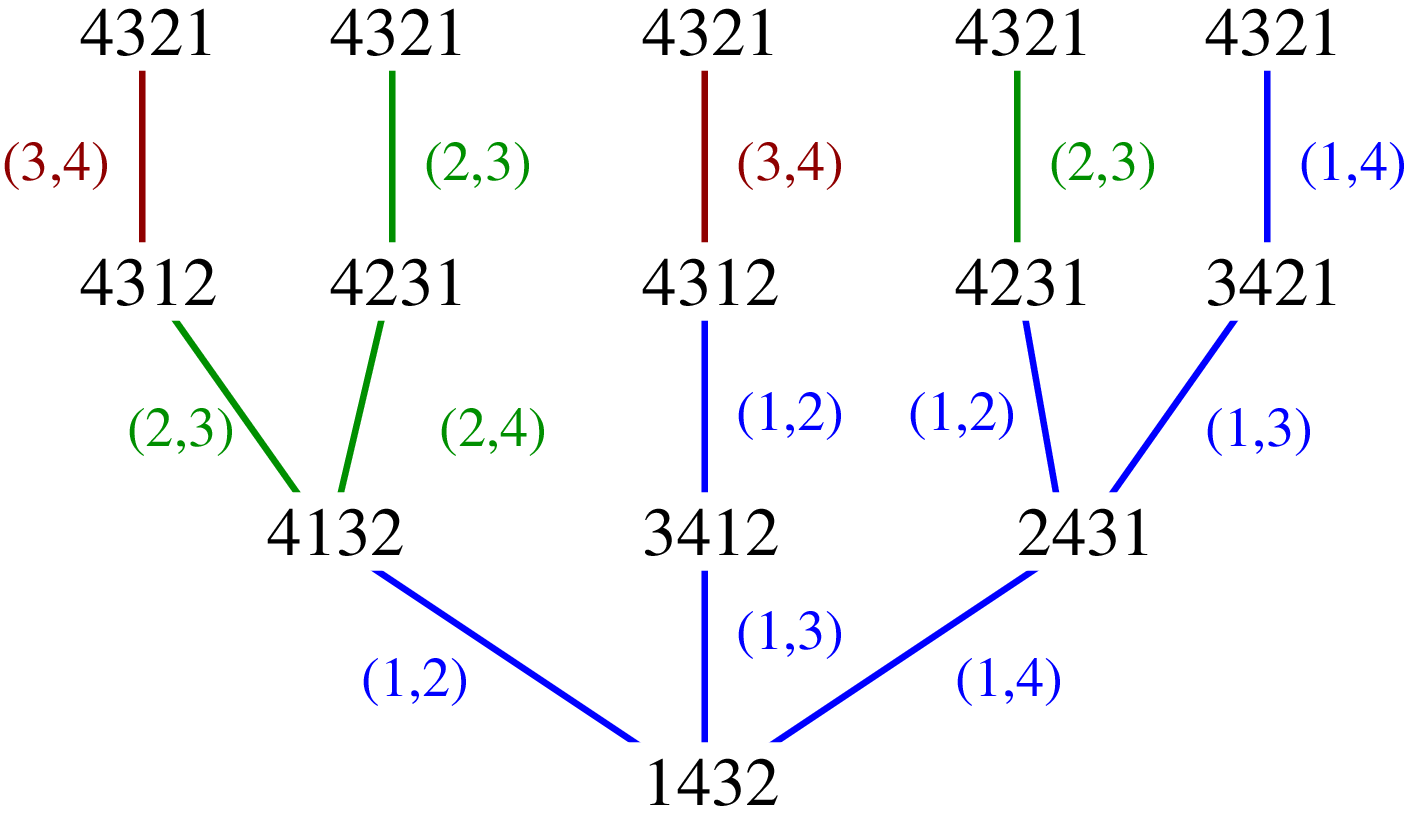}
\]
 Here are the five chains, together with the corresponding terms of
 $\calG_{1432}$.  
\begin{eqnarray*}
 & 1432\Blue{\xrightarrow{\,\underline{(1,2)}\,}}
   4132\ForestGreen{\xrightarrow{\,\underline{(2,3)}\,}}
   4312\Maroon{\xrightarrow{\,\underline{(3,4)}\,}}4321
   & x_1^2x_2\\
 & 1432\Blue{\xrightarrow{\,\underline{(1,2)}\,}}
   4132\ForestGreen{\xrightarrow{\,\underline{(2,4)}\,}}
   4231\ForestGreen{\xrightarrow{\,(2,3)\,}}4321
   & x_1^2x_3 - x_1^2x_2x_3\\
 & 1432\Blue{\xrightarrow{\,\underline{(1,3)}\,}}
   3412\Blue{\xrightarrow{\,(1,2)\,}}
   4312\Maroon{\xrightarrow{\,\underline{(3,4)}\,}}4321
   & x_1x_2^2 - x_1^2x_2^2\\
 & 1432\Blue{\xrightarrow{\,\underline{(1,4)}\,}}
   2431\Blue{\xrightarrow{\,(1,2)\,}}
   4231\ForestGreen{\xrightarrow{\,\underline{(2,3)}\,}}4321
   & x_1x_2x_3 - x_1^2x_2x_3\\
  & 1432\Blue{\xrightarrow{\,\underline{(1,4)}\,}}
   2431\Blue{\xrightarrow{\,(1,3)\,}}
   3421\Blue{\xrightarrow{\,(1,2)\,}}4321
   & x_2^2x_3 - 2x_1x_2^2x_3 + x_1^2x_2^2x_3
\end{eqnarray*}

Consider how each chain may be marked to obtain
a climbing marked chain.
The underlined covers must be marked in any climbing marked
chain.
Every chain can be given three marks, and these have respective
weights $(1,1,1)$, $(1,2,0)$, $(2,0,1)$, $(2,1,0)$, and $(3,0,0)$.
These correspond to the initial monomials with a positive coefficient.
There is a unique way to mark each of the second, third, and fourth chains
twice, and these have weights of $(1,1,0)$, $(1,0,1)$, and $(1,1,0)$,
respectively. 
There are two ways to mark the fifth chain, and each has weight $(2,0,0)$.
These marked chains correspond to the monomials with a negative 
coefficient.
Finally, only the fifth chain admits only a single mark, and this has a weight
of $(1,0,0)$. 
It corresponds to the final monomial in $\calG_{1432}$.

The initial positive monomials and the chains with three
markings correspond to the lowest order terms in the Grothendieck 
polynomial, which is the corresponding Schubert polynomial.
In that case, this example recovers Example~5.4 of~\cite{BS02}.
\end{ex}

We now consider the formula of Lascoux and Sch\"{u}tzenberger 
for decomposing a Grothendieck polynomial $\calG_w(y,x_1,\ldots,x_{n-1})$ into a
linear combination of terms of the form
$y^j\calG_u(x_1,\ldots,x_{n-1})$~\cite{LS82b}, 
using certain operators.
We show that this formula is equivalent to the formula in
Theorem~\ref{T:substitution-formula}.  
This makes explicit Lascoux and Sch\"{u}tzenberger's comment that their formula
is dual to a Pieri-type formula in $K$-theory.

For each simple reflection $s_i$, we define an operator $\phi_i$ on the 
ring $\Lambda[y]$ which commutes with $y$ by describing its action on the
Schubert basis 
\begin{equation}\label{def:phi}
\phi_i \cdot \calG_w = \left\{\begin{array}
                   {c@{\quad  \quad} l}
                   \calG_{s_iw}-\calG_w & \text{ if } s_iw<w, \\
                   0         & \text { otherwise.}\rule{0pt}{15pt}
                   \end{array} \right.
\end{equation}
If $w\in \calS_n$, then let $\mbox{st}(w)\in\calS_{n-1}$ be the permutation 
defined by 
\begin{equation}
{\rm st}(w)(i) = \left\{\begin{array}
        {l@{\quad  \quad} l}
        w(i{+}1)-1& \text{ if } w(i{+}1)>w(1), \\
        w(i{+}1)  & \text { if } w(i{+}1)<w(1)\,.
        \end{array} \right.
\end{equation}
Note that ${\rm st}(w)=w|_{\{2,3,\dotsc,n\}}$.
For example, if $w=365124 \in \calS_6$, then 
$\mbox{st}(w)=54123\in\calS_5$.

\begin{thm}[\cite{LS82b}, Theorem 2.5]\label{T:LS}
 Let $w \in \calS_n$ and set $j=w(1)$.
 Then
 \begin{equation*}
  \calG_w(y,x_1, \ldots, x_{n-1})\ =\
   \prod_{i=j}^{n-2}(1+y\phi_i)\ 
   \cdot y^{j-1}\calG_{{\rm st}(w)}(x_1,\ldots,x_{n-1})\,.
 \end{equation*}
\end{thm}

If we expand the product of Theorem~\ref{T:LS} and apply the operators $\phi_i$,
we get a sum of terms of the form $\pm y^a\calG_v(x)$.
Here is the precise statement.

\begin{lemma}\label{L:sum}
 Let $u\in\calS_{n-1}$.
 Then
\[
   \prod_{i=j}^{n-2}(1+y\phi_i)\ \cdot y^{j-1}\calG_u(x_1,\ldots,x_{n-1})
   \ =\ 
   \sum_{P,Q} (-1)^{|P|} y^{j-1+|P|+|Q|} \calG_{v(P)}(x)\,,
\]
 the sum over all pairs of disjoint subsets $P,Q$ of $\{j,\dotsc,n-2\}$,
 where, if we write $P=\{P_p<P_{p-1}<\dotsb<P_2<P_1\}$ and set 
 $P_0=n-1$ and $P_{p+1}=j-1$, then 
 \begin{equation}\label{E:sum_condition}
  \raisebox{25pt}{\begin{minipage}[t]{350pt}\vspace{-5pt}
    \[
      u\ \gtrdot\ s_{P_1}u\ \gtrdot s_{P_2}s_{P_1}u\ \gtrdot\ \dotsb
       \ \gtrdot\ s_{P_p}\dotsb s_{P_2}s_{P_1}u\ =:\ v(P)\,,
    \]
    is a decreasing chain, and if $q\in Q$ with 
    $P_{i+1}<q<P_i$, then
    \[
      s_q\cdot s_{P_i}\dotsb s_{P_2}s_{P_1}u
      \ \lessdot\ s_{P_i}\dotsb s_{P_2}s_{P_1}u\,.
    \]
  \end{minipage}}
 \end{equation}
\end{lemma}

\begin{proof}
 We first expand the product
\[
 \prod_{i=j}^{n-2}(1+y\phi_i)\ \cdot y^{j-1}\calG_u(x_1,\ldots,x_{n-1})
   \ =\ 
   \sum_{R\subset\{j,j+1,\dotsc,n-2\}}
     y^{j-1+|R|}\phi_{R_r}\dotsb\phi_{R_2}\phi_{R_1} \calG_u(x)\,,
\]
 where we write a subset $R$ as $R_r<\dotsb<R_2<R_1$.
 Each term of the product
 $\phi_{R_r}\dotsb\phi_{R_2}\phi_{R_1} \calG_u(x)$ corresponds to a weakly
 decreasing chain 
\[
  u\ =\ u_0\ \geq\ u_1\ \geq\ \dotsb\ \geq\ u_r\,,
\]
 where, for each $i$, $ u_{i-1}\gtrdot s_{R_i}u_{i-1}$ and 
 $u_i\in\{u_{i-1},s_{R_i} u_{i-1}\}$.
 Let $P\subset R$ correspond to positions of strict descent and set 
 $Q=R\setminus P$.
 Then the final permutation $u_r$ depends only on $P$ and the sign of
 $\calG_{u_r}$ is $(-1)^Q$.
 This completes the proof.
\end{proof}

We define a bijection $\psi$ between the increasing marked chains of
Theorem~\ref{T:substitution-formula} and the pairs $(P,Q)$
satisfying~\eqref{E:sum_condition} of Lemma~\ref{L:sum}. 
Suppose that
 \begin{equation}\label{E:IMC}
  \gamma\ \colon\ w=w_0\ \xrightarrow{\,(1,b_1)\,}\ w_1\ 
     \xrightarrow{\,(1,b_2)\,}\ \dotsb\ 
     \xrightarrow{\,(1,b_t)\,}\ w_t\ =\ n.v
 \end{equation}
is an increasing marked chain in the 1-Bruhat order on $\calS_n$.
Let $j=w(1)$ and for $i=0,1,\dotsc,t$, set $\tau_i:=w_i(1)-1=w(b_i)-1$.
Define
\[
  T\ :=\ \{\tau_1<\tau_2<\dotsb<\tau_t\}\ \subset\ 
          \{j-1,j,j+1,\dotsc,n-2\}\,.
\]
Let $Q\subset T$ consist of those $\tau_i$ for which the cover $(1,b_i)$ is
unmarked in $\gamma$, and set $P:=\{j,j+1,\dotsc,n-2\}\setminus T$.
Set $\psi(\gamma):=(P,Q)$.
If $(P,Q)$ is a valid index for a term in Lemma~\ref{L:sum}, observe that 
\[
  \ell(\gamma)-m(\gamma)\ =\ |Q|,\qquad\mbox{and}\qquad
  n-1-m(\gamma)\ =\ j-1+|P|+|Q|\,,
\]
and thus the power of $y$ in the term of
Theorem~\ref{T:substitution-formula} indexed by $\gamma$ equals the power of $y$
in the term of Lemma~\ref{L:sum} indexed by $(P,Q)$, and likewise their signs
are equal.
Theorem~\ref{T:LS} is a consequence of the following lemma.

\begin{lemma}\label{L:bijection}
  The association $\gamma\mapsto \psi(\gamma)=(P,Q)$ is a bijection between
  increasing marked chains in the $1$-Bruhat order on $\calS_n$ ending in
  the permutation $n.v$ and indices of terms in Lemma~$\ref{L:sum}$
  with $v=v(P)$. 
\end{lemma}

\begin{proof}
 Since the pair $(P,Q)=\psi(\gamma)$ determines the chain $\gamma$,
 we only need to prove that $(P,Q)$ indexes a term in the sum of
 Lemma~\ref{L:sum}, that every such index arises in this way, and that $v=v(P)$.
 The first step is the following consequence of our definitions:

\begin{lemma}
 Suppose that $w\lessdot v$ in the $1$-Bruhat order on $\calS_n$ and that $j=w(1)$
 and $k=v(1)$.
 Let $u={\rm st}(w)\in\calS_{n-1}$.
 Then
\[
   u\ \gtrdot\ s_{k-2}u\ \gtrdot\ \dotsb\ \gtrdot\ 
    s_j s_{j+1}\dotsb s_{k-2} u\ =\ {\rm st}(v)
\]
 is a decreasing chain in $\calS_{n-1}$.
 Note that ${\rm st}(w)={\rm st}(v)$ if $k=j+1$.
\end{lemma}

 Suppose that $\gamma$ is the increasing marked chain~\eqref{E:IMC}.
 Set $j=w(1)$ and $w=w_0$.
 Then
\[
   {\rm st}(w_0)\ \to\ {\rm st}(w_1)\ \to\ \dotsb\ 
      \to\ {\rm st}(w_t)\ =\ v
\]
 is a (weakly) decreasing chain in $\calS_{n-1}$.
 We have
\begin{eqnarray*}
  {\rm st}(w_i)\cdot{\rm st}(w_{i-1})^{-1}&=&
   s_{\tau_{i-1}+1}\cdot s_{\tau_{i-1}+2}\cdot \dotsb\cdot s_{\tau_i-1}\\
  &=& (\tau_{i-1}+1, \tau_{i-1}+2,\dotsb,\tau_i-1)\,,
\end{eqnarray*}
 a cycle.
 Note that the transpositions $s_{\tau_i}$ are exactly those which do
 not occur in any of these cycles.
 Thus, for different $i$, these cycles are disjoint, and we may rearrange them
 to get a decreasing chain in $\calS_{n-1}$ as follows.
 Set $u={\rm st}(w_0)$.
 Then this chain is
\[ 
   u\ \gtrdot\ s_{P_1} u\ \gtrdot\ s_{P_2}s_{P_1}u\ \gtrdot\ \dotsb\ 
        \gtrdot\ s_{P_t}\dotsb s_{P_2}s_{P_1} u=v\,,
\]
where
\[
   \{P_p<\dotsb<P_2<P_1\}\ =\ \{j-1,j,j+1,\dotsc,n-2\}\setminus T
\]
 is the first set in the pair $\psi(\gamma)$.

 Similarly, given a decreasing chain in  $\calS_{n-1}$ from ${\rm st}(w)$ to $v$
 given by left multiplication of adjacent transpositions taken in decreasing
 order (right to left) of their indices (a set $P$), we may invert
 this process to obtain an increasing chain $\gamma$ in the 1-Bruhat order from
 $w$ to $n.v$.
 This shows that the association $\gamma\mapsto P$ is a bijection.
 We complete the proof of the Lemma by addressing the marking/unmarking of
 covers.

 Let $w_{a-1}\xrightarrow{\,(1,b_a)\,}w_a$ be a cover in $\gamma$ and recall
 that $\tau_a=w_a(1)-1$.
 Let $i$ be the index of $P$ such that $P_{i+1}<\tau_a<P_i$.
 We claim that $(1,b_a)$ may be unmarked in $\gamma$ if and only if
 \begin{equation}\label{E:Marking_Condition}
  s_{\tau_a}\cdot s_{P_i}s_{P_{i-1}}\dotsb s_{P_2}s_{P_1}{\rm st}(w)\ \lessdot\ 
  s_{P_i}s_{P_{i-1}}\dotsb s_{P_2}s_{P_1}{\rm st}(w)\,.
 \end{equation}
 This claim will complete the proof of the lemma.
 First, it shows that $\psi(\gamma)=(P,Q)$ is a valid index from
 Lemma~\ref{L:sum}.
 Second, it shows that $\psi$ is surjective.
 Indeed, given a valid index $(P,Q)$, $P$ determines an increasing chain in the
 1-Bruhat order on $\calS_n$  from $w$ to $n.v$.
 According to the claim, if we mark all covers of this chain except those
 $w_{a-1}\xrightarrow{\,(1,b_a)\,}w_a$ where $\tau_a=w_a(1)-1\in Q$, then we
 obtain a valid increasing marked chain $\gamma$ with $\psi(\gamma)=(P,Q)$.

 We prove the claim.
 Since
\[
  \{P_i<P_{i-1}<\dotsb<P_2<P_1\}\ =\ 
  \{\tau_a, \tau_a+1, \dotsc, n-2\}\setminus
  \{\tau_c\mid c\geq a\}\,,
\]
 the leftmost cycle in the product $s_{P_i}\dotsb s_{P_2}s_{P_1}$ of
 disjoint cycles is $(\tau_a+1,\dotsc,\tau_{a+1}-1)$.
 Thus, in $u:=s_{P_i}s_{P_{i-1}}\dotsb s_{P_2}s_{P_1}{\rm st}(w)$, $\tau_a$
 occupies the position that $w_{a+1}(1)-1=w(b_{a+1})-1$ occupies in ${\rm st}(w)$, 
 which is $b_{a+1}-1$.
 Similarly, as $\tau_a<P_i$, it occupies the same positions in $u$ and in 
 ${\rm st}(w)$, namely $b_a-1$.
 Thus
 \begin{eqnarray*}
  (1,b_a)\mbox{ can be unmarked in $\gamma$}
      &\Leftrightarrow&b_a>b_{a+1}\\
      &\Leftrightarrow& \tau_a\mbox{ is to the left of $\tau_a+1$ in $u$}\\
      &\Leftrightarrow& s_{\tau_a} u < u\,,
 \end{eqnarray*}
 which proves the claim.
\end{proof}

\section{Comparison of constructions of Grothendieck polynomials and
   $\calH$-polynomials}\label{S:six} 

We now relate the chain-theoretic construction of Grothendieck polynomials in
Section~\ref{S:five} to some other constructions
and present similar constructions for the $\calH$-polynomials (that were mentioned at the end of Section \ref{S:three}). In fact, in this section we work in a more general context, by referring to the double Grothendieck and $\calH$-polynomials.
We consider the classical construction in~\cite{FK94,LS82b},
reinterpreted by Knutson and Miller~\cite{KM03b}, 
as well a construction of Buch, Kresch, Tamvakis, and
Yong~\cite{BKTY03}. 

Lascoux and Sch\"utzenberger~\cite{LS82b} defined the {\bf double Grothendieck
polynomials}  $\calG_w(x;y)$, which are polynomials in two sets of variables, $x$ 
and $y$, that represent Schubert classes in the $T$-equivariant $K$-theory of
the flag manifold; to be more precise, Lascoux and Sch\"utzenberger used different sets of variables, but a simple change of variables leads to the definition in this paper.
The double Grothendieck polynomials specialize to the ordinary ones upon
setting $y=0$. 
Their definition is nearly identical to that of the ordinary Grothendieck
polynomials, namely
\begin{equation}\label{def:groth}
  \calG_w(x;y)\ :=\ \pi_{w^{-1}\omega_0}\calG_{\omega_0}(x;y)\,,
  \;\;{\rm where}\;\;
   \calG_{\omega_0}(x;y)\ :=\ \prod_{i+j\le n}(x_i+y_j-x_iy_j)\,.
\end{equation}

Fomin and Kirillov~\cite{FK94} use the 0-Hecke algebra
$H_n(0)$ in the study of Grothendieck polynomials. This algebra has generators $u_1,\ldots,u_{n-1}$ satisfying the 
relations
\begin{align*}
  &u_iu_j=u_ju_i\,,\;\;\;\;|i-j|\ge 2\,;\\
  &u_iu_{i+1}u_i=u_{i+1}u_iu_{i+1}\,;\\
  &u_i^2=-u_i\,.
\end{align*}
For $w\in\calS_n$, we can (unambiguously) define $e_w:=u_{i_1}\ldots u_{i_q}$, where
$(i_1,\ldots,i_q)$ is any reduced word for $w$.
The elements $e_w$ form a basis for $H_n(0)$.
Multiplication in $H_n(0)$ is given by 
 \begin{equation}\label{multhecke}
   e_wu_i= \left\{ \begin{array}{ll} e_{ws_i} &\mbox{if}\;\;\;
                         \ell(ws_i)>\ell(w) \\ 
                  -e_{w} &\mbox{otherwise\,,} \end{array} \right.
 \end{equation}
and similarly for left multiplication by $u_i$. 

\begin{remark}{\rm  It is not hard to see that by mapping $u_i$ to the operator $\mu_i=\pi_i-1$ (or $v_i:=u_i+1$ to the operator $\pi_i$), cf. (\ref{def:pi}), we obtain an action of $H_n(0)$ on $\Z[x_1, \ldots, x_n]$. Furthermore, by mapping $u_i$ to the operator $\phi_i$ defined in (\ref{def:phi}), we obtain an action of $H_n(0)$ on $\Lambda:=\Z[x_1, \ldots, x_n]/ I_n\simeq H^*(\Fl V)$. }
\end{remark}

Consider the following expression
 \begin{equation}\label{E:Cauchy}
  \calG(x;y):=\prod_{i=1}^{n-1}\prod_{j=n-1}^{i}\big(1+(x_i+y_j-x_iy_j)u_j\big)\,.
 \end{equation}
The order of the indices ($i$ increases and $j$ decreases) is
important, as $H_n(0)$ is noncommutative. 
Fomin and Kirillov showed that 
 \begin{equation}\label{gfgroth}
   \calG(x;y)=\sum_{w\in\calS_n}\calG_w(x;y)e_w\,.
 \end{equation}

Double Schubert polynomials have a similar expression~\cite{FK96,FS95},
replacing $x_i+y_j-x_iy_j$ by $x_i-y_j$ and $H_n(0)$ by the nilCoxeter
algebra. 
That formula implies their construction in terms of reduced words
and compatible sequences in~\cite{BJS93,LS82b},
which is represented graphically in terms of certain line diagrams for
permutations introduced in~\cite{FK96} and called {\bf resolved braid configurations}. The same objects were called {\bf rc-graphs} in~\cite{BeBi93}.  
We repeat this for double Grothendieck and double
$\calH$-polynomials.
Each of the formulas we give below implies that the lowest
homogeneous components of these polynomials are the double Schubert
polynomials.  

Rc-graphs are certain subsets of
 \[
   \{(i,j)\in\Z_{>0}\times\Z_{>0}\mid i+j\le n\}\,
 \]
associated to permutations. 
Any such subset of pairs $R$ is linearly ordered
 \begin{equation}\label{ordcross}
    (i,j)\le (i',j') \ \Longleftrightarrow\ (i< i')
    \quad\mbox{or}\quad (i=i'\ \mbox{and}\  j\ge j')\,.
 \end{equation}
Let $(i_k,j_k)$ be the $k$th pair in this linear order, and let
$a_k:=i_k+j_k-1$. 
Consider the sequences (words)
\[
  {\rm word}(R):=(a_1,a_2,\ldots)\quad{\rm and}\quad {\rm  comp}(R):=(i_1,i_2,\ldots)\,.
\]
An rc-graph $R$ is associated to $w$ if ${\rm word}(R)$ is a reduced
word for $w$.
When this happens, ${\rm comp}(R)$ is called a {\bf compatible
sequence} with ${\rm word}(R)$. 
We write $\calR(w)$ for the rc-graphs associated to a
permutation $w$, and $w(R)$ for the permutation associated to a given rc-graph $R$. 

We represent an rc-graph $R$ graphically as a planar history of the
inversions of $w$ (see Example~\ref{exrcg}). 
To this end, we draw $n$ lines going up and to the right, such that
the $i$th line starts at position $(i,1)$ and ends at position
$(1,w_i)$.
Two lines meeting at position $(i,j)$ cross if and only if $(i,j)$
is in $R$.

\begin{ex}\label{exrcg}
 Here are two rc-graphs associated to the permutation $215463$.  

\[
   \begin{picture}(300,120)
   \put(2,1){\epsfxsize=3.8in \epsffile{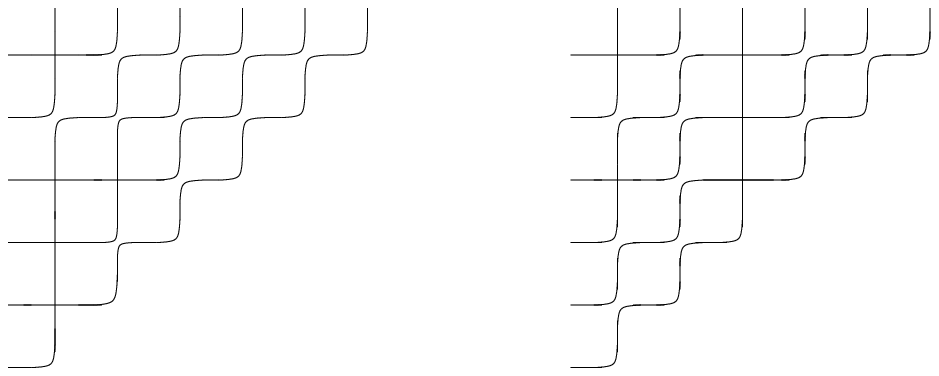}}
   \put(20,111){1}\put(38,111){2}\put( 56,111){3}
   \put(73,111){4}\put(92,111){5}\put(110,111){6}
   \put(182,111){1}\put(200,111){2}\put(218,111){3}
   \put(235,111){4}\put(254,111){5}\put(272,111){6}
   \put(0,90){1}   \put(162,90){1}
   \put(0,72){2}   \put(162,72){2}
   \put(0,54){3}   \put(162,54){3}
   \put(0,36){4}   \put(162,36){4}
   \put(0,18){5}   \put(162,18){5}
   \put(0,0){6}    \put(162,0){6}
   \end{picture}
\]
\end{ex}

We use {\bf right} and {\bf left marked rc-graphs} to construct double Grothendieck
polynomials. 
These objects are equivalent to the {\bf unreduced pipe dreams} in~\cite{KM03b}.
A right marked rc-graph $R\,\hat\,$ is an ordered bipartition 
$R\,\hat\,=(R,R^*)$ where $R$ is an
rc-graph, and $R^*$ is a set (possibly empty) of marked positions
$(i,j)\not\in R$ (as before, $i+j\le n$) such that the two lines avoiding each
other at $(i,j)$ cross in another position, northeast of $(i,j)$. 
We call such positions {\bf right absorbable positions}; note that they were
called {\bf absorbable elbow tiles} in~\cite{KM03b}. We write $(i,j)\in R\,\hat\,=(R,R^*)$ for $(i,j)\in R\sqcup R^*$. 
Let ${\rm rabs}(R)$ be the collection of such positions for an 
rc-graph $R$.
Then $R^*\subseteq{\rm rabs}(R)$. 
As noted in~\cite{KM03b}, there is a unique rc-graph in 
$\calR(w)$ with ${\rm rabs}(R)=\emptyset$, the
{\bf bottom rc-graph} of~\cite{BeBi93}.
The rc-graph on the left in Example~\ref{exrcg} is a bottom rc-graph. 
Let $\calR\hat\,(w)$ be set of right marked rc-graphs with
underlying rc-graph $R$ in $\calR(w)$. 
The positions in a right marked rc-graph $R\,\hat\,$ are ordered as
in~\eqref{ordcross}, and we can thus define the word of $R\,\hat\,$, which is denoted by ${\rm word}(R\,\hat\,)$. By right multiplying the generators of $H_n(0)$ corresponding to the word of $R\,\hat\,=(R,R^*)$ - based on~\eqref{multhecke}, the length increases by 1 precisely when the corresponding position in $R\,\hat\,$ belongs to $R$; in fact, this justifies the term ``right marked rc-graph''. We conclude that $R\,\hat\,\in\calR\hat\,(w)$ if and
only if ${\rm word}(R\,\hat\,)$ evaluates to $\pm e_w$ in $H_n(0)$. 
We combine this with~\eqref{E:Cauchy} and~\eqref{gfgroth} (by expanding
the product and collecting terms) and obtain the following construction.

\begin{thm}\cite{FK94,KM03b}\label{markedrc} 
 We have
 \begin{equation}\label{expg1}
     \calG_w(x;y)=\sum_{R\,\hat\,\in\calR\hat\,(w)}
       (-1)^{|R\,\hat\,|-\ell(w)}
       \prod_{(i,j)\in R\,\hat\,}(x_i+y_j-x_iy_j)\,.
 \end{equation}
 On the other hand, we also have
 \begin{equation}\label{expg2}
     \calG_w(x;y)=\sum_{R\in\calR(w)}
              \prod_{(i,j)\in R}(x_i+y_j-x_iy_j)
            \prod_{(i,j)\in{\rm rabs}(R)}(1-x_i)(1-y_j)\,.
 \end{equation}
\end{thm}

The second formula is a consequence of the first: summing over all
right marked rc-graphs with the same underlying rc-graph $R$ means
summing over all subsets of ${\rm rabs}(R)$. 
Note that, in fact, it is (\ref{gfgroth}) which was stated in \cite{FK94} (for
ordinary Grothendieck polynomials); but, as we discussed above and is also
briefly mentioned in \cite{FK94}, (\ref{gfgroth}) immediately leads to the
formulas in Proposition \ref{markedrc}. Later, these formulas appeared in
\cite{KM03b} in terms of unreduced pipe dreams; the proofs in \cite{KM03b} were
based on the simplicial topology of subword complexes.   

\begin{rem}\label{R:Left}
 The notions defined above have left analogs, which are based on left
 multiplication in $H_n(0)$. 
 For instance, a {\bf left marked rc-graph} $\hat\,R$ is
 an ordered bipartition $\hat\,R=(R,\,^*\!R)$, where $R$ is
 an rc-graph, and $^*\!R$ is a set (possibly empty) of marked
 positions $(i,j)\not\in R$ (as before, $i+j\le n$) such that the two lines
 avoiding each other at  $(i,j)$ cross in another position, 
 {\it southwest} of $(i,j)$. 
 We likewise have left absorbable positions and the set 
 $\hat\:\calR(w)$ of left marked rc-graphs underlying some
 rc-graph in $\calR(w)$. 
 We are lead to a nearly identical construction of the double
 Grothendieck polynomials, cf. Example \ref{exg1}. Note that each subset $\widehat{R}$ of $\{(i,j)\in\Z_{>0}\times\Z_{>0}\mid i+j\le n\}$ can be viewed as a left marked rc-graph $\hat\,R=(R,\,^*\!R)$ and as a right marked rc-graph $R\,\hat\,=(R',R^*)$ such that $\widehat{R}=R\sqcup\,^*\!R=R'\sqcup R^*$. Since formula (\ref{expg1}), which is in terms of right marked rc-graphs $R\,\hat\,=(R,R^*)$,  only depends on $R\sqcup R^*$, its analog for left marked rc-graphs is identical. 
\end{rem}

\begin{ex}\label{exg1} 
  The  permutation $132$ has two rc-graphs.
  We draw them below; their right and left absorbable positions are
  marked by circles and squares, respectively.  

\[
    \epsfxsize=1.7in \epsffile{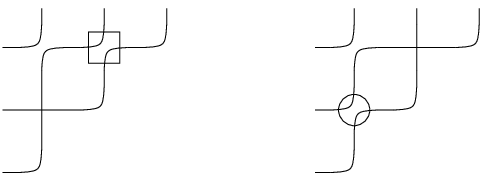}
\]
 This gives two formulas for $\calG_{132}(x;y)$, based on right and left
 absorbable positions, respectively: 
 \begin{align*}
  \calG_{132}(x;y)\ &=\ (x_2+y_1-x_2y_1)+(x_1+y_2-x_1y_2)(1-x_2)(1-y_1)\\
  &=\ (1-x_1)(1-y_2)(x_2+y_1-x_2y_1)+(x_1+y_2-x_1y_2)\,.
 \end{align*}
\end{ex}

\begin{cor}\label{gsymm} We have
 \[
  \calG_w(x;y)=\calG_{w^{-1}}(y;x)\,.
 \]
\end{cor}
\begin{proof} 
  Transposing a right marked rc-graph gives a bijection between 
 $\calR\hat\,(w)$ and $\hat\:\calR(w^{-1})$. The corollary follows by combining this bijection with~\eqref{expg1}.
\end{proof}

We now deduce a related formula for double Grothendieck polynomials
in terms of {\bf double rc-graphs}. 
A double rc-graph~\cite{BeBi93} for a permutation $w$ is a pair
$(u,v)$ of permutations such that $v^{-1} u=w$ and
$\ell(w)=\ell(u)+\ell(v)$, together with an rc-graph for each of $u$
and $v$. 
We write the rc-graph for $v$ upside down and above that of $u$,
and label its rows with negative numbers (see Example~\ref{exg2}). 
Write $\calD(w)$ for the set of all double rc-graphs $D$ for
$w$. 
Right and left marked double rc-graphs and right and left absorbable
positions in double rc-graphs are defined as above, and we extend
the entire system of notation.

\begin{thm}[Corollary 3 of~\cite{BKTY03}, see Remark~\ref{bktycons}]\label{bkty} 
   Set $x_i:=y_{-i}$ for $i<0$. We have
 \begin{equation}\label{dexpg1}
     \calG_w(x;y)=
    \sum_{D\,\hat\,\in\calD\hat\,(w)}(-1)^{|D\,\hat\,|-\ell(w)}
     \prod_{(i,j)\in D\,\hat\,}x_i\,.
 \end{equation}
 There is a similar formula for left marked double rc-graphs.
 We also have
 \begin{equation}\label{dexpg2}
   \calG_w(x;y)=
   \sum_{D\in\calD(w)}\prod_{(i,j)\in D}x_i
   \prod_{(i,j)\in{\rm rabs}(D)}(1-x_i)\,.
 \end{equation}
 There is a similar formula, obtained by replacing ${\rm rabs}(D)$ by ${\rm labs}(D)$.
\end{thm}

\begin{proof}
These follow from the Cauchy identity for Grothendieck polynomials~\cite{FK94}:
 \[
   \calG_w(x;y)\ =\ 
   \sum_{e_v\cdot e_u=\pm e_w}(-1)^{\ell(u)+\ell(v)-\ell(w)}\calG_u(x;0)\,\calG_v(0;y)\,,
 \]
where the product $e_v\cdot e_u$ is taken in the 0-Hecke
algebra $H_n(0)$. 
More precisely, this also uses~\eqref{expg1} with
$y_j=0$ and $\calG_v(0;y)=\calG_{v^{-1}}(y;0)$ (Corollary~\ref{gsymm}).
\end{proof}

\begin{rem}\label{bktycons}
  In a slightly different form, formula~\eqref{dexpg1} is
  Corollary~3 of~\cite{BKTY03}, where it was obtained by specializing a
  $K$-theoretic quiver formula of Buch~\cite{Bu02}.  
\end{rem}

\begin{ex}\label{exg2}
  There are four double rc-graphs for the permutation $132$
  which we indicate below, marking their right absorbable positions
  by circles. 
\[
   \epsfxsize=3.8in 
  \epsffile{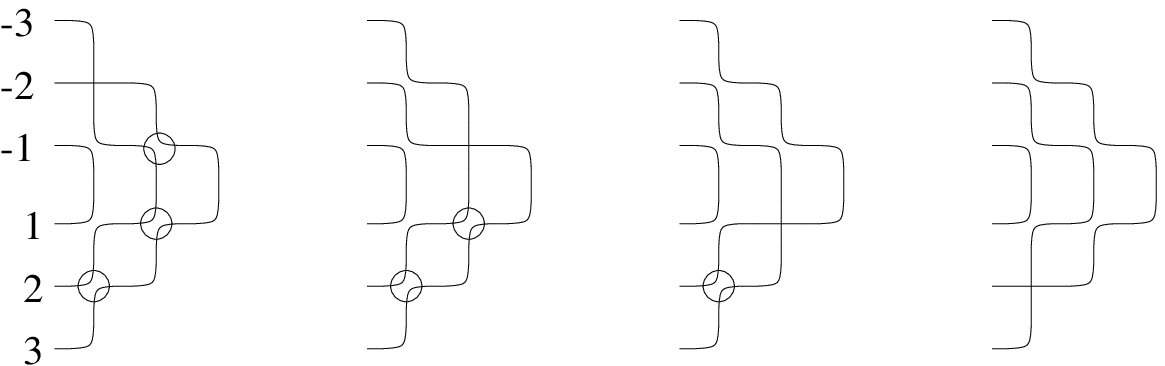}
\]
Hence, we have the following formula for $\calG_{132}(x;y)$:
\[y_2(1-y_1)(1-x_1)(1-x_2)+y_1(1-x_1)(1-x_2)+x_1(1-x_2)+x_2\,.\]
\end{ex}

We relate the construction of double Grothendieck polynomials in
Theorem~\ref{markedrc} to the chain-theoretic construction in
Section~\ref{S:five}. 
This uses a bijection given in~\cite{LS03} between rc-graphs
and chains in Bruhat order satisfying a certain condition.
Our conventions are different than in~\cite{LS03}, so we must prove that this
bijection has the desired properties.  

Let $R$ be an rc-graph in $\calR(w)$ and define a saturated chain 
$\gamma(R)$ in Bruhat order from $w$ to $\omega_0$ as follows.
First, construct a sequence of rc-graphs 
$R=R_0,R_1,\ldots,R_{\ell(\omega_0)-\ell(w)}$: given $R_i$, add the
pair $(k,a)$ to $R_i$ where  $(k,a)$ is the pair not in $R_i$ that
is minimal in the lexicographic order on pairs. 
By Lemma~10 of~\cite{LS03}, 
$\gamma(R)\colon w=w(R_0)\lessdot w(R_1)\lessdot\dotsb\lessdot
w(R_{\ell(\omega_0)-\ell(w)})=\omega_0$ is a saturated chain in
Bruhat order with covers labeled by the pairs $(k,a)$ chosen in its construction.
The lexicographic minimality of the label covers imply that $\gamma(R)$
is a climbing chain.  
Now, let $\hat\,R=(R,\,^*\!R)\in\hat\:\calR(w)$.
Set the climbing marked chain $\gamma(\hat\,R)$ to be the chain $\gamma(R)$
where covers $(k,a)$ not in $^*\!R$ are marked.

\begin{prop}\label{bijection} 
 The association $\gamma$ is a bijection between 
 $\hat\:\calR(w)$ and the set of climbing marked chains 
 from $w$ to $\omega_0$. 
\end{prop}

\begin{proof}
 We first check that the map $R\mapsto\gamma(R)$ is a bijection between
 $\calR(w)$ and climbing chains from $w$ to $\omega_0$. 
 Given a climbing chain 
 $\gamma\colon w=w_0\lessdot w_1\lessdot \dotsb\lessdot w_p=\omega_0$, 
 set $t_{k_i,m_i}:=w_{i-1}^{-1}w_i$ and $t_{a_i,c_i}:=w_iw_{i-1}^{-1}$, with
 $k_i<m_i$ and $a_i<c_i$. 
 It suffices to show that the pairs $(l,b)$ with $l+b\le n$ that are not
 labels $(k_i,m_i)$ of covers in $\gamma$ form an rc-graph for $w$. 
 We induct on $p=\ell(\omega_0)-\ell(w)$ with the case $p\le 1$
 immediate. 
 If $p>1$, consider the rc-graph $R_1$ with $w(R_1)=w_1$ corresponding to the
 climbing chain $w_1\lessdot  \dotsb\lessdot w_p=\omega_0$. 
 By induction, this chain is $\gamma(R_1)$ and $R_1$ contains all positions
 $(l,b)$ lexicographically smaller than $(k_2,a_2)$. 
 Since $\gamma$ was a climbing chain,  either $k_1=k_2$ and $a_1<a_2=c_1$, or
 else $k_1<k_2$ and $c_1=n+1-k_1$. 
 Therefore, position $(k_1,a_1)$ is in $R_1$. 
 Removing the crossing at $(k_1,a_1)$ gives a line diagram for the permutation
 $w$ (since $w_1(m_1)=a_1$).
 Also all positions directly above and to the left of $(k_1,a_1)$ lie in this
 diagram, as shown on the left below. 
 Since this line diagram has $\ell(w)=\ell(w_1)-1$ crosses, it must be an rc-graph. 
\[
  \begin{picture}(300,135)
   \put(1,0){\epsfxsize=4.1in \epsffile{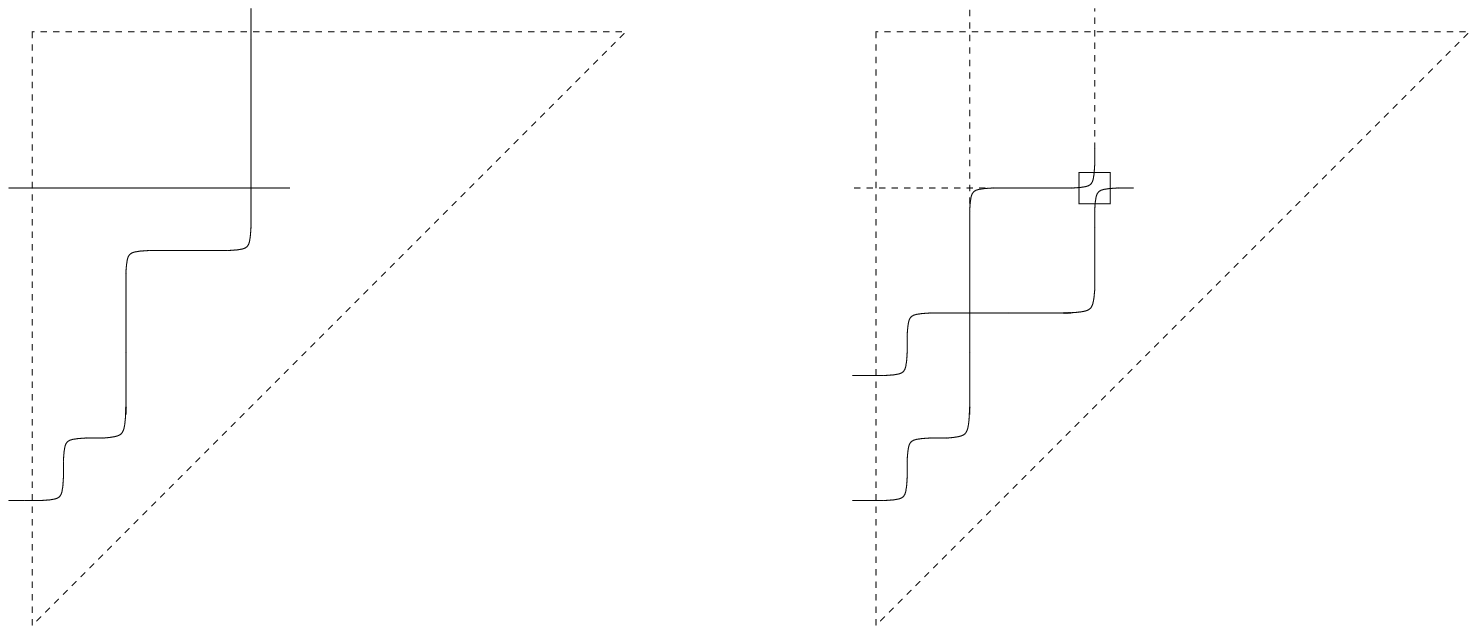}}
   \put(0,24){$m_1$} \put(3,82){$k_1$} \put(58,125){$a_1$}
   \put(162,24){$m_2$} \put(168,82){$k$}
   \put(162,47){$m_1$} \put(198,125){$a'$} \put(222,125){$a$}
  \end{picture}
\]

Consider a left marked rc-graph $\hat\,R=(R, \,^*\!R)$ with 
associated climbing marked chain $\gamma(R)$. 
Let $(k,a)$ be a pair in $\,^*\!R$ so that the two lines avoiding each other at
that position cross in a position southwest of $(k,a)$. 
Let $m_1<m_2$ be the two rows where these lines end (on the
left side of $R$, see the second figure above). Consider the
rightmost position $(k,a')$ in row $k$ which does not belong to $R$
and is to the left of $(k,a)$; clearly, such a position exists. 
The positions $(k,a')$ and $(k,a)$ correspond to successive covers
labeled $(k,m_2)$ and $(k,m_1)$ in the associated chain; the second
cover is not marked, but this portion of the chain is increasing in
the sense of Section~\ref{S:five}. 
On the other hand, every cover which is not marked arises in this
way.  This concludes the proof.  
\end{proof}

The bijection of Proposition~\ref{bijection} combined with
Theorem~\ref{markedrc} gives a formula for double 
Grothendieck polynomials in terms of climbing marked chains. 
Given a cover $u\lessdot v$ in a chain,
set $t_{i,k}=u^{-1}v$ and $t_{j,l}=vu^{-1}$, with $i<k$ and $j<l$. 
Given a climbing marked chain $\gamma$, let 
\[\mu_{x,y}(\gamma):=\prod(x_i+y_j-x_iy_j)\,,\]
the product is taken over all marked covers in $\gamma$.

\begin{cor}\label{groth2chains}
 We have
 \begin{equation}\label{exch}
   \calG_w(x;y)\ =\ \sum_{\gamma}(-1)^{\ell(\gamma)-m(\gamma)}\,
          \calG_{\omega_0}(x;y)/\mu_{x,y}(\gamma)\,,
 \end{equation}
 the sum over all climbing marked chains from $w$ to $\omega_0$. 
 Setting $y_j=0$ recovers the formula of Theorem~$\ref{T:monomial-formula}$.  
\end{cor}

Let us now turn to the double $\calH$-polynomials. These were defined (for each $w\in\calS_n$) in~\cite{La90b} in a completely similar way to the double Grothendieck polynomials (but using different variables than the ones we use). To be more precise, one only needs to replace the operators $\pi_{w^{-1}\omega_0}$ with $\mu_{w^{-1}\omega_0}$ in (\ref{def:groth}); these operators were defined at the end of Section \ref{S:three} for a polynomial ring in variables $x_i$, and the definition is identical when there is a second set of variables $y_j$. Unlike the double Grothendieck polynomials, which are stable under the inclusion $\calS_n\hookrightarrow\calS_{n+1}$, the double $\calH$-polynomials are not. Setting $y_j=0$ gives the ordinary $\calH$-polynomials, which have been previously considered at the end of Section \ref{S:three}. 

In order to study the double $\calH$-polynomials, we need to consider new generators $v_i:=u_i+1$ of the 0-Hecke algebra $H_n(0)$. These satisfy the same
relations as $u_i$, except that $v_i^2=v_i$ instead of $u_i^2=-u_i$. 
Similarly to the definition of the basis $\{{e}_w\mid w\in\calS_n\}$ of $H_n(0)$, the generators $v_i$ determine a basis $\{\widetilde{e}_w\mid w\in\calS_n\}$. More precisely, given $w\in\calS_n$ and any reduced word $(i_1,\ldots,i_q)$ for $w$, we can (unambiguously) define $\widetilde{e}_w:=v_{i_1}\ldots v_{i_q}$. Clearly, there is an algebra involution on $H_n(0)$ given by $u_i\mapsto -v_i$. The formulas in the following lemma express the corresponding change of bases. This lemma also appears as Lemma 1.13 in \cite{La90b}, where it is proved based on induction and the Exchange Lemma for Coxeter groups. We present here a direct proof; we believe that it offers more insight, once the appropriate definitions are made. 

\begin{lemma}\cite{La90b}\label{change-base} We have
\[
  {e}_w\ =\ \sum_{w'\le w}(-1)^{\ell(w)-\ell(w')}\widetilde{e}_{w'}\,,\;\;\;\widetilde{e}_w\ =\ \sum_{w'\le w}{e}_{w'}\,.
 \]
\end{lemma}

\begin{proof}
Fix a reduced word for $w$, say $Q=(i_1,\ldots,i_q)$, and 
recall some notions from~\cite{KM03b} related to subword complexes. 
Given a subword $P\subseteq Q$, say
$P=(i_{j_1},\ldots,i_{j_p})=(k_1,\ldots,k_p)$, let
$\delta(P):=v_{k_1}\dotsb v_{k_p}$ in $H_n(0)$. If $P$ is reduced,
we say that $i_r\in Q\setminus P$ is {\bf absorbable} if the subword
$T:=P\cup i_r$ of $Q$ satisfies 
$\delta(T)=\delta(P)$, and if the unique $i_t$ satisfying
$\delta(T\setminus i_t)=\delta(P)$ has $t<r$ (Lemma 3.5.2 in~\cite{KM03b}
guarantees existence and uniqueness). 
Given $w'\le w$, a reduced word for $w'$ appears as a subword of $Q$;
conversely, any subword of $Q$ evaluates to some $\widetilde{e}_{w'}$ in
$H_n(0)$, where $w'\le w$. 
Fix $w'\le w$ for the remainder of the proof. 
The {\bf facet adjacency graph}~\cite{KM03b} is a directed graph whose vertices
are reduced words for $w'$ that are subwords of $Q$.
It has a directed edge from $P'$ to $P$ whenever 
$P'=(P\cup i_r)\setminus i_t$, where $i_r$ is absorbable and $i_t$ is as
above. 
This graph is acyclic~\cite[Remark 4.5]{KM03b}, with a unique source. 

The coefficient of $\widetilde{e}_{w'}$ in 
\[
  e_w\ =\ (v_{i_1}-1)\dotsb(v_{i_q}-1)\ =\ 
  \sum_{T\subset Q} (-1)^{q-|T|}\delta(T)\,,
\]
is the sum of $(-1)^{q-|T|}$ for each $T\subseteq Q$ with
$\delta(T)=\widetilde{e}_{w'}$.
Consider the sum over those $T$ having the same lexicographically first reduced
subword $P$ for $w'$ will be 0 whenever $P$ has at least one absorbable element
in $Q$. 
Indeed, such $T$ are obtained by adjoining to $P$ a subset of its absorbable
elements. 
This completes the proof of the first formula, as there is a unique such $P$ with no absorbable elements
(the unique source of the facet adjacency graph). The second formula follows by applying the involution given by $u_i\mapsto -v_i$ to the first formula.
\end{proof}

The first formula in Lemma \ref{change-base} immediately leads to a formula for the double $\calH$-polynomials in terms of the double Grothendieck polynomials, which generalizes (\ref{def:hpoly}); recall that the latter formula is the combinatorial counterpart of the geometric statement in Proposition \ref{P:duality}. Furthermore, by substituting the first formula in Lemma \ref{change-base} into~\eqref{gfgroth}, we can realize the generating function $\calG(x;y)$ for the double Grothendieck polynomials as a generating function for the double $\calH$-polynomials. We present these results in the corollary below.

\begin{cor}\label{cor:hpoly} We have
\[
   \calH_w(x;y)\ :=\ \sum_{\omega_n\geq v\ge w}(-1)^{\ell(v)-\ell(w)}\calG_v(x;y)\,.
 \]
We also have
\[
   \calG(x;y)\ =\ \sum_{w\in\calS_n}\calH_w(x;y)\widetilde{e}_w\,.
 \]
 \end{cor}

Let us rewrite $\calG(x;y)$ as
 \begin{align*}\calG(x;y)\ &:=\ 
     \prod_{i=1}^{n-1}
      \prod_{j=n-1}^i[1+(x_i+y_j-x_iy_j)(v_j-1)]\\
    \;&=\ \prod_{i=1}^{n-1}
          \prod_{j=n-1}^i[(1-x_i)(1-y_j)+(x_i+y_j-x_iy_j)v_j]\,.
 \end{align*}
The second formula in Corollary \ref{cor:hpoly} now leads to a formula for $\calH_w(x;y)$ in terms of right
(respectively left) marked rc-graphs which is similar to~\eqref{expg1}; more precisely, we have
 \begin{equation}\label{hmarkedrc}
  \calH_w(x;y)\ =\ \sum_{R\,\hat\,\in\calR\hat\,(w)}
      \prod_{(i,j)\in R\,\hat\,}(x_i+y_j-x_iy_j)
      \prod_{(i,j)\not\in R\,\hat\,,\,i+j\le n}(1-x_i)(1-y_j)\,.
 \end{equation}
Once again, summing over all right marked rc-graphs with the same underlying
rc-graph $R$ means summing over all subsets of ${\rm rabs}(R)$. 
But we have
 \[
   \sum_{A\subseteq{\rm rabs}(R)}
   \prod_{(i,j)\in A}(x_i+y_j-x_iy_j)
   \prod_{(i,j)\in{\rm rabs}(R)\setminus A}(1-x_i)(1-y_j)\ =\ 1\,.
\]
Hence we obtain the following formula for the $\calH$-polynomials.

\begin{cor}\label{exph} We have 
\begin{equation}\label{exph16}
 \calH_w(x;y)\ =\ \sum_{R\in\calR(w)}\prod_{(i,j)\in R}(x_i+y_j-x_iy_j)
 \prod_{(i,j)\not\in R\sqcup{\rm rabs}(R),\,i+j\le n}(1-x_i)(1-y_j)\,.
 \end{equation}
There is a similar formula, obtained by replacing ${\rm rabs}(R)$ by ${\rm
  labs}(R)$.  
\end{cor}

\begin{ex}\label{exh1}
  Recall from Example~\ref{exg1} the two rc-graphs for the permutation $132$,
  including their right and left absorbable positions. We have the following two
  formulas for $\calH_{132}(x;y)$, based on right and left absorbable positions,
  respectively: 
\begin{align*}
  &(1-x_1)(1-y_1)\cdot\big((1-x_1)(1-y_2)(x_2+y_1-x_2y_1)+(x_1+y_2-x_1y_2)\big)\\
  &=(1-x_1)(1-y_1)\cdot\big((x_2+y_1-x_2y_1)+(1-x_2)(1-y_1)(x_1+y_2-x_1y_2)\big)\,.
\end{align*}
\end{ex}

Combining the formula for $\calH_w(x;y)$ in Corollary~\ref{exph}
with the formula~\eqref{expg2} for $\calG_w(x;y)$, we obtain the
following relationship between the two polynomials. Note that this relationship appears in a different guise in \cite{La90b}, as Lemma 2.9 2). 

\begin{cor}\label{C:614}\cite{La90b}
We have
\[
  \calH_w(x;y)\ =\ 
  (-1)^{\ell(w)}\calH_{(1,\ldots,n)}(x;y)\,\calG_w(\widetilde{x};\widetilde{y})\,,
\]
where $\widetilde{x}_i=x_i/(x_i-1)$ and $\widetilde{y}_i=y_i/(y_i-1)$.
\end{cor}

Using Corollary~\ref{C:614}, we show that the double
$\calH$-polynomials have similar properties to double Grothendieck
polynomials. 

\begin{cor}\label{prophpoly}
We have $\calH_w(x;y)=\calH_{w^{-1}}(y;x)$. 
We also have the Cauchy formula
\[\calH_w(x;y)\ =\ \sum_{v\cdot u=w}\calH_u(x;0)\,\calH_v(0;y)\,.\]
\end{cor}

Using~\eqref{hmarkedrc} and Corollary~\ref{prophpoly}, we obtain the analog
of~\eqref{dexpg2} for double $\calH$-polynomials.

\begin{cor}
 Setting $x_i:=y_{-i}$ for $i<0$, we have
 \begin{equation}\label{dexph16}
  \calH_w(x;y)\ =\ \sum_{D\in\calD(w)}
      \prod_{(i,j)\in D}x_i
      \prod_{(i,j)\not\in D\sqcup{\rm rabs}(D)}(1-x_i)\,.
 \end{equation}
There is a similar formula obtained by replacing ${\rm rabs}(D)$ by ${\rm labs}(D)$.
\end{cor}

\begin{ex} 
 Recall the double rc-graphs for the permutation $132$ from Example~\ref{exg2}.
 These give the following formula for $\calH_{132}(x;y)$
\[
  (1-y_1)(1-x_1)\cdot\big(y_2+y_1(1-y_2)+x_1(1-y_2)(1-y_1)
  +x_2(1-y_2)(1-y_1)(1-x_1)\big)\,.
\]
\end{ex}

We deduce from~\eqref{exph16} and the bijection of Proposition~\ref{bijection} a
formula for the double $\calH$-polynomials in terms of climbing chains. 
This requires a different notion of a marked chain.
The covers in a climbing chain will now be marked so that when we cut each
subchain in a given $k$-Bruhat order at the permutations before each marking, the
subchain decomposes into \emph{decreasing} segments 
(see Section~\ref{S:five} for the definition of and increasing/decreasing chain). 
By contrast, the original notion of marking required the resulting segments to
be increasing and it also required the first cover in a chain to be marked. 
This new marking procedure can be rephrased as follows: a cover labeled $(k,j)$
has to be marked if it is preceded by one labeled $(k,i)$ with $i>j$.

\begin{cor}
 We have
 \begin{equation}\label{exhch}
  \calH_w(x;y)\ =\ \sum_{\gamma}(-1)^{\ell(\gamma)-m(\gamma)}\,
   \calH_{\omega_0}(x;y)/\mu_{x,y}(\gamma)\,,
 \end{equation}
 the sum over all climbing chains from $w$ to $\omega_0$ with this new marking.
\end{cor}

\begin{ex}
 Recall the five climbing chains in Example~\ref{fivechains} corresponding to
 the permutation ${1432}$. 
 We write the terms of $\calH_{1432}(x)=\calH_{1432}(x;0)$ corresponding to
 each. 
 The underlined covers are those that must be marked. 
\begin{eqnarray*}
 & 1432\Blue{\xrightarrow{\,(1,2)\,}}
   4132\ForestGreen{\xrightarrow{\,(2,3)\,}}
   4312\Maroon{\xrightarrow{\,(3,4)\,}}4321
   & x_1^2x_2(1-x_1)(1-x_2)(1-x_3)\\
 & 1432\Blue{\xrightarrow{\,(1,2)\,}}
   4132\ForestGreen{\xrightarrow{\,(2,4)\,}}
   4231\ForestGreen{\xrightarrow{\,\underline{(2,3)}\,}}4321
   & x_1^2x_3(1-x_1)(1-x_2)\\
 & 1432\Blue{\xrightarrow{\,(1,3)\,}}
   3412\Blue{\xrightarrow{\,\underline{(1,2)}\,}}
   4312\Maroon{\xrightarrow{\,(3,4)\,}}4321
   & x_1x_2^2(1-x_1)(1-x_3)\\
 & 1432\Blue{\xrightarrow{\,(1,4)\,}}
   2431\Blue{\xrightarrow{\,\underline{(1,2)}\,}}
   4231\ForestGreen{\xrightarrow{\,(2,3)\,}}4321
   & x_1x_2x_3(1-x_1)(1-x_2)\\
  & 1432\Blue{\xrightarrow{\,(1,4)\,}}
   2431\Blue{\xrightarrow{\,\underline{(1,3)}\,}}
   3421\Blue{\xrightarrow{\,\underline{(1,2)}\,}}4321
   & x_2^2x_3(1-x_1)
\end{eqnarray*}
\end{ex}

There is one construction of double Grothendieck polynomials which we do not
know how to relate to the classical construction in terms of rc-graphs.
This is due to Lascoux and involves certain alternating sign matrices~\cite{La02}.
Our attempts to do so suggest that any bijective relation between these formulas
will be quite subtle.

\section{Specialization of Grothendieck polynomials and permutation 
patterns}\label{S:seven_geom}

Until Remark~\ref{R:identities}, we fix subsets $P$ and $Q$ of $[m{+}n]$ such
that $P\sqcup Q=[m{+}n]$ with $\#P=m$ and $\#Q=n$.
We will compute the effect of the map $\varphi_{P,Q}$ on the Schubert basis of
the Grothendieck rings and use those results to obtain formulas for the
specialization of Grothendieck polynomials at two sets of variables.

In the situation of Section~\ref{S:four}, the Billey-Postnikov pattern
map~\cite{BiPo} is given by  
 \begin{equation}\label{E:pattern}
   \calS_{m+n}\ni w\ \longmapsto\ (w|_P, w|_Q)\in\calS_m\times\calS_n\,,
 \end{equation}
where $w|_P$ is the permutation in $\calS_m$ with the same set of 
inversions as
the sequence $(w_{p_1}, w_{p_2}, \dotsc, w_{p_m})$.
The fibers of this surjective map are identified with 
decompositions of $[m{+}n]$ into disjoint subsets
of sizes $m$ and $n$.
The identification is through the values taken by $w$ at the positions
of $P$ and $Q$.

\begin{dfn}[Section 4.5 of~\cite{BS98}]\label{D:varepsilon}
 Define permutations $\varepsilon(P,Q)$ and $\zeta(P,Q)$ in $\calS_{m+n}$ by
 requiring that both have image $(e,e)\in\calS_m\times\calS_n$ under the pattern
 map~\eqref{E:pattern}, and the values they take at $P$ and
 at $Q$ to be as follows
 \[
  \begin{array}{rclcl}
     \varepsilon(P,Q)&\colon&P\longmapsto [m]&\mbox{ and }&
                      Q\longmapsto \{m{+}1\dotsc,m{+}n\}\\
     \zeta(P,Q)     &\colon&Q\longmapsto [n]&\mbox{ and }&\rule{0pt}{13pt}
                      P\longmapsto \{n{+}1\dotsc,m{+}n\}\\
  \end{array}
 \]
 Suppose that $m=4$ and $n=3$ with $P=\{1,2,4,7\}$ and $Q=\{3,5,6\}$.
 Then
\[
 \varepsilon(P,Q)\ =\  \Blue{12}\Magenta{5}\Blue{3}\Magenta{67}\Blue{4}
  \qquad\mbox{and}\qquad
 \zeta(P,Q)\ =\  \Blue{45}\Magenta{1}\Blue{6}\Magenta{23}\Blue{7}\,.
\]

 For $u\in\calS_m$ and $v\in\calS_n$, define
 $v\times u := (v_1,\dotsc,v_n,n+u_1,\dotsc,n+u_m) \in\calS_{m+n}$.
 Then we have two sections of the pattern map defined by
\[
  \begin{aligned}
   \calS_m\times\calS_n\ \ni\  (u,v)&\longmapsto&
     u\times v\cdot\varepsilon(P,Q)\ \in\ \calS_{m+n}\\
   \calS_m\times\calS_n\ \ni\  (u,v)&\longmapsto& 
    v\times u\cdot\zeta(P,Q)\ \in\ \calS_{m+n}\\
  \end{aligned}
\]
Let $u=\Blue{2413}$ and $v=\Magenta{312}$.
Then 
$u\times v=\Blue{2413}\Magenta{756}$ and 
$v\times u=\Magenta{312}\Blue{5746}$, and with $P$ and $Q$ as above, we have
\[
 u\times v\cdot \varepsilon(P,Q)\ =\ 
   \Blue{24}\Magenta{7}\Blue{1}\Magenta{56}\Blue{3}
 \qquad\mbox{and}\qquad
 v\times u \cdot \zeta(P,Q)\ =\ 
  \Blue{57}\Magenta{3}\Blue{4}\Magenta{12}\Blue{6}\,.
\]

\end{dfn}

Some of the following elementary properties of these permutations may be found
in Section 4.5 of~\cite{BS98}, and the rest we leave to the reader.

\begin{lemma}\label{L:interval}
 We have
\begin{enumerate}
 \item[(a)] For any $u\in\calS_m$ and $v\in\calS_n$, we have 
     \[
        u\times v\cdot\varepsilon(P,Q)\ =\  
       \omega_{m+n}\cdot (\omega_nv\times \omega_mu)\cdot\zeta(P,Q)\,.
     \]
 \item[(b)] The map
    \[
      \calS_m\times\calS_n\ \ni\ (u,v)\ 
      \longmapsto\ v\times u\cdot\zeta(P,Q)
    \]
     induces an order preserving isomorphism between 
     $\calS_m\times\calS_n$ and the interval 
    \[ 
       [\zeta(P,Q),\ \omega_n\times\omega_m\cdot\zeta(P,Q)]
     \]
      in the Bruhat order.
 \item[(c)] 
      In particular, if $\alpha\in\calS_m$ and $\beta\in\calS_n$, then 
     \[
       \ell(\beta\times\alpha\cdot\zeta(P,Q))\: -\:   
         \ell(v\times u\cdot\zeta(P,Q))\ =\ 
         \ell(\alpha)\:+\: \ell(\beta)\:-\:\bigl(\ell(u)\:+\: \ell(v)\bigr)\,.
      \]
\end{enumerate}
\end{lemma}

\begin{prop}[Corollary~4.5.2 of~\cite{BS98}]\label{P:BS}
 Let $\Edot,\Epdot\subset U$ and $\Fdot,\Fpdot\subset V$ be pairs of
 opposite flags.
 If we set $R:=\{m+1,\dotsc,m+n\}$ and 
 $S:=\{n+1,\dotsc,m+n\}$, then 
\[
  \Gdot\ :=\ \varphi_{S,[n]}(\Edot,\Fdot)\qquad\mbox{and}\qquad
  \Gpdot\ :=\ \varphi_{[m],R}(\Epdot,\Fpdot)\,,
\]
 are opposite flags in $U\oplus V$.
 Furthermore, for any $u\in\calS_m$ and $v\in\calS_m$, 
 we have
 \begin{equation}\label{E:PM}
  \varphi_{P,Q}(X_u\Edot\times X_v\Fdot)\ =\ 
     X_{v\times u\cdot\zeta(P,Q)}(\Gdot)
    \ \cap \ 
     X_{\varepsilon(P,Q)}(\Gpdot)\,.
 \end{equation}
\end{prop}

In our runing example, this is
\[
  \varphi_{P,Q}(X_{\Blue{2413}}\Edot\times X_{\Magenta{312}}\Fdot)\ =\ 
     X_{\Blue{57}\Magenta{3}\Blue{4}\Magenta{12}\Blue{6}}(\Gdot)
    \ \cap \ 
     X_{\Blue{12}\Magenta{5}\Blue{3}\Magenta{67}\Blue{4}}(\Gpdot)\,.
\]

Proposition~\ref{P:BS} follows from Lemma 4.5.1 of~\cite{BS98}, which shows that
the left side of~\eqref{E:PM} is contained in the right side by comparing the 
linear-algebraic definitions  of the Schubert varieties and the definitions of
$\varepsilon(P,Q)$ and $\zeta(P,Q)$.
Equality follows as both sides are irreducible and have the
same dimension. 
This dimension calculation is essentially Lemma~\ref{L:interval}(c).

We deduce the effect of the pushforward $(\varphi_{P,Q})_*$ on the
Schubert basis.
We use the notation of Section~\ref{S:four}.
By the K\"unneth theorem,
\[
  K^0(\calF')\ =\ 
  K^0(\Fl U \times \Fl V) \ \simeq\ 
  K^0(\Fl U)\otimes K^0(\Fl V)\,.
\]
Here, the isomorphism is induced by the identifications of sheaves
\[
  \calO_{X_u\times X_v}\ \simeq\ \pi_1^*\calO_{X_u}\otimes_{\calO_{\calF'}}
   \pi_2^*\calO_{X_v}\,,
\]
where $\pi_1,\pi_2$ are the projections onto the factors of 
$\calF'=\Fl U \times \Fl V$.
On Schubert classes this is 
$[\calO_{X_u\times X_v}]\mapsto[\calO_{X_u}]\otimes[\calO_{X_v}]$.
Thus 
\[ 
   \{[\calO_{X_u}]\otimes[\calO_{X_v}]\mid u\in\calS_m, \ 
      v\in\calS_n\}
\]
is a basis for $K^0(\Fl U \times \Fl V)$.
Proposition~\ref{P:BS} has the immediate corollary.

\begin{cor}~\label{C:pushForward}
 For $u\in\calS_m$ and $v\in\calS_m$,  we have
\[
  (\varphi_{P,Q})_*\bigl( [\calO_{X_u}]\otimes[\calO_{X_v}]\bigr)\ =\ 
  [\calO_{X_{\varepsilon(P,Q)}\cap X'_{v\times u\cdot\zeta(P,Q)}}]\ =\ 
  [\calO_{X_{\varepsilon(P,Q)}}]\cdot[\calO_{X_{v\times u\cdot\zeta(P,Q)}}]\,.
\]
\end{cor}

We also compute the pushforward of the ideal sheaves
$[\calI_u]\otimes[\calI_v]$.

\begin{thm}\label{T:pushForward}
 For $u\in\calS_m$ and $v\in\calS_n$, 
 \begin{equation}\label{E:pushForward}
  (\varphi_{P,Q})_*\bigl([\calI_u]\otimes[\calI_v]\bigr)\ =\
   [\calI_{v\times u\cdot \zeta(P,Q)}]\cdot [\calO_{X_{\varepsilon(P,Q)}}]\,. 
 \end{equation}
\end{thm}

In our running example, 
\[
  (\varphi_{P,Q})_*
    \bigl([\calI_{\Blue{2413}}]\otimes[\calI_{\Magenta{312}}]\bigr)\ =\
   [\calI_{\Blue{57}\Magenta{3}\Blue{4}\Magenta{12}\Blue{6}}]
   \cdot [\calO_{X_{\Blue{12}\Magenta{5}\Blue{3}\Magenta{67}\Blue{4}}}]\,. 
\]

\begin{proof}
 We expand the right side of~\eqref{E:pushForward} using
 Proposition~\ref{P:duality}. 
 Set $\zeta=v\times u\cdot\zeta(P,Q)$ and 
 $\varepsilon=\varepsilon(P,Q)$.  
 Then 
\begin{eqnarray}
  [\calI_\zeta]\cdot[\calO_{X_\varepsilon}]
   &=& \sum_{w\geq\zeta} (-1)^{\ell(w\zeta)}
            [\calO_{X_w}]\cdot[\calO_{X_\varepsilon}]\nonumber\\
   &=& \sum_{\omega_{m+n}\cdot\varepsilon\geq w\geq\zeta}
              (-1)^{\ell(w\zeta)}
            [\calO_{X_w}]\cdot[\calO_{X_\varepsilon}]\,,
            \label{EA:sum}
\end{eqnarray}
 as $[\calO_{X_w}]\cdot[\calO_{X_\varepsilon}]=0$ unless 
 $\omega_{m+n}\cdot\varepsilon\geq w$.
 Using Lemma~\ref{L:interval}, which characterizes the interval between 
 $\zeta=v\times u\cdot\zeta(P,Q)$ and 
 $\omega_n\times\omega_m\cdot\zeta(P,Q)=\omega_{m+n}\cdot\varepsilon(P,Q)$,
 and in particular the assertion concerning lengths in part (c) (modulo
 2), the sum~\eqref{EA:sum} becomes
 \begin{multline*}
  \quad \sum_{\omega_n\geq\beta\geq v} \ \sum_{\omega_m\geq\alpha\geq u}
     (-1)^{\ell(\beta v)} (-1)^{\ell(\alpha u)} 
    \bigl[\calO_{X_{\beta\times\alpha\cdot\zeta(P,Q)}}\bigr]\cdot
      [\calO_{X_{\varepsilon(P,Q)}}]\\
  \begin{aligned}
    &\ =\   \sum_{\beta\geq v}\sum_{\alpha\geq u} 
     (-1)^{\ell(\beta v)} (-1)^{\ell(\alpha u)}
  (\varphi_{P,Q})_*\bigl([\calO_{X_\alpha}]\otimes[\calO_{X_\beta}]\bigr)
     \quad\\
    &\ =\   (\varphi_{P,Q})_*
      \Bigl(\sum_{\alpha\geq u} (-1)^{\ell(\alpha u)} 
          [\calO_{X_\alpha}]\Bigr)
      \otimes
      \Bigl(\sum_{\beta\geq v} (-1)^{\ell(\beta v)} 
          [\calO_{X_\beta}]\Bigr)\ ,\quad
   \end{aligned}
 \end{multline*}
 which equals $(\varphi_{P,Q})_*\bigl([\calI_u]\otimes[\calI_v]\bigr)$,
 by Proposition~\ref{P:duality} again.
\end{proof}

We compute the pullbacks of $[\calO_{X_w}]$ and $[\calI_w]$ along the map 
$\varphi_{P,Q}$.

\begin{thm}\label{T:decomposition_coefficients}
 For $w\in\calS_{m+n}$, we have
\begin{eqnarray*}
  (\varphi_{P,Q})^*[\calO_{X_w}]&=& \sum_{(u,v)\in\calS_m\times\calS_n}
           c^{u\times v\cdot \varepsilon(P,Q)}_{w,\,\varepsilon(P,Q)}\,
             [\calO_{X_u}]\otimes[\calO_{X_v}]\\
  (\varphi_{P,Q})^*[\calI_w]&=&  \sum_{(u,v)\in\calS_m\times\calS_n}
           c^{\omega_{m+n} w}_{\varepsilon(P,Q),\ 
           (\omega_nv\times\omega_mu)\cdot\zeta(P,Q)} \,
             [\calI_u]\otimes[\calI_v] \\
\end{eqnarray*}
\end{thm}

\begin{proof}
  Since the classes $[\calO_{X_u}]\otimes[\calO_{X_v}]$ form a basis for 
  $K^0(\Fl U \times \Fl V)$, there are integers $a^{u,v}_w$
  defined by the identity
\[
  (\varphi_{P,Q})^*[\calO_{X_w}]\ =\  \sum_{(u,v)\in\calS_m\times\calS_n}
           a^{u, v}_w\, [\calO_{X_u}]\otimes[\calO_{X_v}]\,.
\]

 We extract the coefficient $a^{u,v}_w$ from this
 expression using the dual basis $[\calI_u]\otimes[\calI_v]$ for
 $K^0(\Fl U\times \Fl V)$ and map $\chi$.
 Thus 
\begin{eqnarray*}
  a^{u,v}_w&=&\chi\bigl( \varphi^*_{P,Q}([\calO_{X_w}])\cdot
   [\calI_{\omega_mu}]\otimes[\calI_{\omega_n v}]\bigr)\\
  &=& \chi\bigl([\calO_{X_w}]\cdot(\varphi_{P,Q})_*
       \bigl([\calI_{\omega_mu}]\otimes
              [\calI_{\omega_nv}]\bigr)\bigr)\,,
\end{eqnarray*}
 by the projection formula in $K$-theory~\cite[5.3.12]{CG97}.
 By Theorem~\ref{T:pushForward}, we obtain
 \begin{equation}\label{E:constant}
  a^{u,v}_w\ =\ \chi\bigl([\calO_{X_w}]\cdot
       [\calI_{(\omega_nv\times \omega_mu)\cdot\zeta(P,Q)}]
       \cdot[\calO_{X_{\varepsilon(P,Q)}}]\bigr)\,.
 \end{equation}
 As $\omega_{m+n}^2=e$, the identity permutation, Lemma~\ref{L:interval}(a)
 implies that
 \begin{equation}\label{E:identification}
  (\omega_nv\times \omega_mu)\cdot \zeta(P,Q)\ =\ 
   \omega_{m+n}\cdot(u\times v)\cdot\varepsilon(P,Q),
 \end{equation}
 and so by the Equation~\eqref{E:StrConst} we
 conclude that 
\[
  a^{u,v}_w\ =\ c^{u\times v\cdot\varepsilon(P,Q)}_{w,\,\varepsilon(P,Q)}\,.
\]

 The second formula is obtained in the same way as the first, except
 that the formula~\eqref{E:constant} for the coefficients is
\[
  \chi\bigl([\calI_w]\cdot
       [\calO_{X_{(\omega_nv\times \omega_mu)\cdot\zeta(P,Q)}}]
       \cdot[\calO_{X_{\varepsilon(P,Q)}}]\bigr)\,.
\]
 As $w=\omega_{m+n}(\omega_{m+n} w)$, this is
 $c^{\omega_{m+n} w}_{\varepsilon(P,Q),\ 
      (\omega_nv\times \omega_mu)\cdot\zeta(P,Q)}$,
 which completes the proof.
\end{proof}

As described in Section~\ref{S:three}, computations in the Grothendieck ring 
yield results about Grothendieck polynomials, by the stability of
Grothendieck polynomials.

\begin{thm}\label{T:SubsComp}
 Let $w\in\calS_{m+n}$ and suppose that $a\geq m$ and $b\geq n$ are integers
 such that  
\[
  \psi_{P,Q}(\calG_w(x))\ \in\ 
   R_a(y)\otimes R_b(z)\,.
\]
 Then, for any $P'\supset P$ and $Q'\supset Q$ with 
 $P'\sqcup Q'=[a+b]$ where $\#P'=a$ and $\#Q'=b$, 
 we have
\[
    \psi_{P,Q}(\calG_w(x))\ =\ 
   \sum_{u\in \calS_a, \ v\in\calS_b}\;
    c^{u\times v\cdot\varepsilon(P',Q')}_{w,\,\varepsilon(P',Q')}\,
    \calG_u(y)\cdot\calG_v(z)\,.
\]

 Similarly, if $a$ and $b$  satisfy
 $\psi_{P,Q}R_m(x)\subset R_a(y)\otimes R_b(z)$, then we have
\[
    \psi_{P,Q}(\calH_w(x))\ =\ 
   \sum_{u\in\calS_a, v\in \calS_b}\;
    c^{\omega_{a+b} w}_{\varepsilon(P',Q'),\ 
       (\omega_bv\times\omega_au)\cdot\zeta(P',Q')}\,
   \calH_u(y)\cdot \calH_v(z)\,.
\]
\end{thm}

\begin{proof}
 Since $w\in\calS_{m+n}$, we have $\calG_w\in\mathbb{Z}[x_1,\dotsc,x_{m+n}]$,
 and so for any $P'\supset P$  and $Q'\supset Q$ we have 
 $\psi_{P,Q}(\calG_w(x))=\psi_{P',Q'}(\calG_w(x))$.
 Proposition~\ref{P:pattern_map} asserts that 
 $\varphi^*_{P',Q'}$ agrees with the action of $\psi_{P',Q'}$ on the
 quotient $\mathbb{Z}[y]/\calI_a\otimes\mathbb{Z}[z]/\calI_b$.
 Since $\psi_{P,Q}(\calG_w(x))$ lies in $R_a(y)\otimes
 R_b(z)$, which is identified with this quotient, the formula for
 $\psi_{P,Q}(\calG_w(x))$ then follows from 
 Theorem~\ref{T:decomposition_coefficients}. 

 The formula for $\psi_{P,Q}(\calH_w(x))$ follows in a similar fashion from 
 Theorem~\ref{T:decomposition_coefficients}.
\end{proof}

\begin{rem}\label{R:identities}
 Different choices of the sets $P'$ and $Q'$ (and of the numbers $a$ and $b$)
 will give different formulas for $\psi_{P,Q}(\calG_w)$.
 These different choices lead to a plethora of identities among the Schubert
 structure constants. 
 We will describe one class
 of such formulas, beginning with a simple form of the statement of
 Theorem~\ref{T:SubsComp}. 

 Note that $a=b=m+n$ will suffice in the hypothesis of
 Theorem~\ref{T:SubsComp}, so that $a{+}b{=}2(m{+}n)$. 
 If we suppress the sizes of the sets $P$ and $Q$, then we may replace $m{+}n$
 by a single number $k$, and we have the following simpler statement.

\begin{prop}
  Suppose that $P,Q\subset[k]$ with $P\sqcup Q=[k]$ and $w\in\calS_k$.
  Then, for any $P'\supset P$ and $Q'\supset Q$ of size $k$ with  
  $P'\sqcup Q'=[2k]$, we have the formula,  
\[
   \psi_{P,Q}(\calG_w(x))\ =\ \sum_{u,v\in\calS_k}
    c^{u\times v\cdot\varepsilon(P',Q')}_{w,\,\varepsilon(P',Q')}\,
    \calG_u(y)\cdot\calG_v(z)\,.
\]
\end{prop}

This gives many identities among the structure constants.
Similar identities for Schubert polynomials were the
motivation in~\cite{BS98} for studying the maps $\varphi_{P,Q}$.

\begin{cor}
 Suppose that $P,P', Q,Q'\subset[2k]$ are sets of size $k$
 such that $P\sqcup Q=P'\sqcup Q'=[2k]$ with $P\cap[k]=P'\cap[k]$ and
 $Q\cap[k]=Q'\cap[k]$.
 Then, for any $u,v,w\in\calS_k$, we have the following identity of Schubert
 structure constants
\[
       c^{u\times v\cdot\varepsilon(P,Q)}_{w,\,\varepsilon(P,Q)} 
  \ =\ c^{u\times v\cdot\varepsilon(P',Q')}_{w,\, \varepsilon(P',Q')}\,. 
\]
\end{cor}

\end{rem}

\begin{ex}\label{Ex:First_substitution}
 We illustrate Theorem~\ref{T:SubsComp} with the substitution
 $\calG_{1432}(z_1,y_1,z_2,z_3)$. 
 That is, with $P=\{2\}$ and $Q=\{1,3,4\}$.
 Since $\psi_{P,Q}(\calG_{1432})\in R_3(y)\otimes R_3(z)$, we may take 
 $P'=\{2,5,6\}$ and $Q'=\{1,3,4\}$.
 Then $\varepsilon(P',Q')=415623$ and we may calculate the product
 $\calG_{143256}\cdot\calG_{415623}$ in the quotient ring
 $\Z[x_1,\dotsc,x_6]/\calI_6$ to obtain
\[
 \begin{aligned}
  &\hspace{1.2em}\Blue{ \calG_{436512}+\calG_{526413}+\calG_{624513}+\calG_{615423}
  -\calG_{536412}-\calG_{634512}-2\calG_{625413}+\calG_{635412}}\\
  &+\calG_{456213}+\calG_{462513}+\calG_{465123}+\calG_{561423}\ 
   -\calG_{456312}-\calG_{463512}-2\calG_{465213}-\calG_{546213} \\
  &-2\calG_{562413}-\calG_{564123}-\calG_{642513}-\calG_{645123}
   -\calG_{651423}\\
  &+\calG_{465312}+\calG_{546312}+\calG_{563412}+2\calG_{564213}
   +\calG_{643512}+2\calG_{645213}+2\calG_{652413}+\calG_{654123}\\
  &-\calG_{564312}-\calG_{645312}-\calG_{653412}-2\calG_{654213}
    +\calG_{654312}\\
 \end{aligned}
\]
 Only the 8 terms in the first row have indices of the form 
 $u\times v\cdot\varepsilon(P',Q')$, and so only those contribute to the expression 
 for $\psi_{P,Q}(\calG_{1432})$ of Theorem~\ref{T:SubsComp}.
 Inspecting this first row, we see that this expression is
\[
  \begin{aligned}
   &\hspace{1.2em}\ 
            \calG_{312}(y)\cdot\calG_{132}(z)
       \ +\ \calG_{213}(y)\cdot\calG_{231}(z)
       \ +\ \hspace{9pt}\calG_{213}(y)\cdot\calG_{312}(z)
       \ +\ \calG_{123}(y)\cdot\calG_{321}(z)\\
   &
       -\ \calG_{312}(y)\cdot\calG_{231}(z)
       \ -\ \calG_{312}(y)\cdot\calG_{312}(z)
 \ -\ 2\,\calG_{213}(y)\cdot\calG_{321}(z)
       \ +\ \calG_{312}(y)\cdot\calG_{321}(z)\,.
 \end{aligned}
\]
 Using the identification of Grothendieck polynomials for $\calS_3$ of
 Example~\ref{E:S3}, we obtain
\[
  \begin{aligned}
   &\ \hspace{1.2em}
    y_1^2\cdot(z_1+z_2-z_1z_2)\ +\ y_1\cdot z_1z_2\ 
 +\ y_1\cdot z_1^2\ +\  1\cdot z_1^2z_2\\
&
   -\ y_1^2\cdot z_1z_2\ -\ y_1^2\cdot z_1^2\ -\ 2 y_1\cdot z_1^2z_2
    \ +\ y_1^2\cdot z_1^2z_2\,,
 \end{aligned}
\]
 Comparing with~\eqref{E:G1432}, we recognize this as 
 $\calG_{1432}(z_1,y_1,z_2,z_3)$. 
\end{ex}

\section{Substitution of a single variable}\label{S:seven}

We give the formula for the substitution of a single variable used in
Section~\ref{S:five}. 
For $u\in\calS_n$ and $j,q\in[n{-}1]$, let 
$\xi_{q,j}(u)\in\calS_\infty$ be the permutation whose first $n{+}1$ values are  
$u_1,\dotsc,u_{q-1},\,n{+}j{+}1$, $u_q,\dotsc,u_n$, and whose remaining values 
are increasing.

\begin{thm}\label{T:SingleSubstitution}
 Let $w\in\calS_n$ and $1\leq q<n$.
 Then
\[
  \calG_w(x_1,\dotsc,x_{q-1},\,y,\,x_{q},\dotsc,x_{n-1})\ =\ 
   \sum_{\substack{j\geqslant 0\\u\in\calS_n}}
    c^{\xi_{q,j}(u)}_{w,r_{[q,n-q+1]}}\, y^j\, \calG_u(x)\,,
\]
 
 If $q=1$, then this simplifies to become
\[
  \calG_w(y, x_1,\dotsc,x_{n-1})\ =\ 
   \sum_{\substack{j\geqslant 0\\u\in\calS_{n-1}}}  
   c^{n{+}j.u}_{w,r_{[1,n-1]}}\, y^j\, \calG_u(x)\,.
\]
\end{thm}

\begin{proof}
 For $w\in\calS_n$, the maximum degree of the $q$th variable in the Grothendieck  
 polynomial $\calG_w(x)$ is $n{-}q$.
 If we define 
\[
   P\ :=\ \{1,2,\dotsc,q{-}1,\, q{+}1,\dotsc,n{+}1\}
      \qquad\textrm{and}\qquad
   Q\ :=\ \{q,\ n{+}2, n{+}3, \dotsc, 2n{-}q{+}1\}\,,
\]
 then 
\[
  \psi_{P,Q}(\calG_w)\ =\ 
  \calG_w(x_1,\dotsc,x_{q-1},\,y,\,x_{q},\dotsc,x_{n-1})\ 
   \in\  R_n(x)\otimes R_{n+1-q}(y)\,.
\]
 By Theorem~\ref{T:SubsComp}, we have 
\begin{equation}\label{E:onevar}
  \psi_{P,Q}(\calG_w)\ =\ 
  \sum_{\substack{u\in\calS_n\\v\in\calS_{n+1-q}}}
   c^{u\times v\cdot\varepsilon(P,Q)}_{w,\,\varepsilon(P,Q)}\,
    \calG_u(x)\cdot \calG_v(y)\,.
\end{equation}

 We complete the proof of the theorem by identifying the permutation
 $u\times v\cdot\varepsilon(P,Q)$.
  By Definition~\ref{D:varepsilon}, 
\[
  \varepsilon(P,Q)\ =\ [1,\dotsc,q{-}1,n{+}1,q,\dotsc,n,n{+}2,\dotsc]\ =\ 
   r_{[q,n-q+1]}\,.
\] 
 Similarly, the permutation $u\times v\cdot\varepsilon(P,Q)$ is 
\[
  u\times v\cdot\varepsilon(P,Q)\ =\ 
   [u_1,\dotsc,u_{q-1},\,n+v_1,\,u_q,\dotsc, u_n, n+v_2,\dotsc,n+v_{n+1-q}]\,.
\]

 Since $\varepsilon(P,Q)=r_{[q,n-q+1]}$, the constants in~\eqref{E:onevar} are
 $c^{u\times v\cdot\varepsilon(P,Q)}_{w, r_{[q,n-q+1]}}$.
 These vanish unless $w<_q u\times v\cdot\varepsilon(P,Q)$.
 Since $w\in\calS_n$, it has no descents after position $n$.
 Then the characterization~\cite[Theorem A]{BS98} of the $q$-Bruhat order
 implies that $u\times v\cdot\varepsilon(P,Q)$ has no descents after position $n$.
 Thus $v$ equals 
 $[j{+}1,1,\dotsc,j,j{+}2,\dotsc,n{-}q{+}1]=r_{[1,j]}$, for some $j$.
 By Example~\ref{E:S3}, $\calG_v(y)=y^j$.
 Also $u\times v\cdot\varepsilon(P,Q)$ is the
 permutation with no descents after position $n$ whose restriction to
 $P$ is $u$, and whose value at position $q$ is $n{+}j{+}1$, for some $j$. 
 That is, $u\times v\cdot\varepsilon(P,Q)=\xi_{q,j}(u)$.
 The first formula follows from this.
 
 If $q=1$, then 
 $\calG_w(y,x_1,\dots,x_{n-1})\in R_{n-1}(x)\otimes R_{n}(y)$,  
 and so we may replace $n$ by $n{-}1$ when referring to the $x$-variables, $u$,
 $P$, etc.
 This gives the second formula.
\end{proof}

\begin{ex}\label{E:SingleSubstitution}
 We illustrate Theorem~\ref{T:SingleSubstitution} in the case when $q=2$
 and $w=1432$. 
 Then $r[2,3]=15234$, and we compute 
 $\calG_{1423}\cdot\calG_{15234}$ to be
\[
  \calG_{1732456} + \calG_{2631457} + \calG_{3612457} + \calG_{3521467}
 -\calG_{2731456} - \calG_{3712456} - 2\,\calG_{3621457}
 +\calG_{3721456}\,.
\]
 Each index in this sum has the form $\xi_{2,j}(u)$ for
 some $u\in\calS_4$, and so each term contributes to the expression
 for $\calG_{1432}(x_1,y_1,x_2,x_3,\dotsc)$.
 The values $(j,u)$ are, from left to right
\[
  (2,1324),\ (1,2314),\ (1,3124),\ (0,3214),\ 
  (2,2314),\ (2,3124),\ (1,3214),\ (2,3214)\,.
\]
 Furthermore, each of these permutations lies in $\calS_3$.
 (This is because $\phi_{P,Q}(\calG_{1432})\in R_3(x)\otimes R_{3}(y)$, rather
 than $R_4(x)\otimes R_{3}(y)$.)
 Thus Theorem~\ref{T:SingleSubstitution} implies the formula for 
 $G_{1432}(x_1,y_1,x_2,x_3,\dotsc)$:
 \[
 \begin{aligned}
   &\hspace{1.2em}\ 
        y^2 \calG_{132}(x)
   \ +\ y   \calG_{231}(x)
   \ +\ \hspace{9pt}
        y   \calG_{312}(x)
   \ +\ 1   \calG_{321}(x)\\
  &  -\ y^2 \calG_{231}(x)
   \ -\ y^2 \calG_{312}(x)
   \ -\ 2y\,\calG_{321}(x)
   \ +\ y^2 \calG_{321}(x)\,,
 \end{aligned}
 \]
 which we recognize from Example~\ref{Ex:First_substitution} as 
 $G_{1432}(x_1,y,x_2,x_3)$.
\end{ex}

\section{General substitution formula}\label{S:eight}

The formula of Theorem~\ref{T:SubsComp} may be iterated to give formulas for
specializing the variables in a Grothendieck polynomial at more than two sets
of variables, expressing such a specialization as a sum of products of Grothendieck
polynomials in each set of variables.
This procedure shows that the coefficients in such an expression are sums of
products of Schubert structure constants.
We use the full geometry of the pattern map to show that 
these coefficients are naturally Schubert structure
constants.

Let $\Adot=(A_1,A_2,\dotsc,A_s)$ be a collection of disjoint subsets of $[m]$
with union $[m]$, that is, a set composition of $[m]$. 
Set $|\Adot|:=(|A_1|, |A_2|, \dotsc, |A_s|)$, a composition of $m$. 
Let ${\bf Y}^{(1)}$, ${\bf Y}^{(2)}$, \dots, ${\bf Y}^{(s)}$ be infinite sets of 
variables.
The specialization map 
 \begin{equation}\label{E:specialization}
   \psi_{\Adot}\ \colon\  
    \mathbb{Z}[x_1,x_2,\dotsc, x_m]\ \longrightarrow\  
    \bigotimes_{i=1}^s \mathbb{Z}[{\bf Y}^{(i)}]\,,
 \end{equation}
is defined by specializing the $x$-variables in the positions of $A_i$ to the
variables in ${\bf Y}^{(i)}$:
\[
  \psi_{\Adot}\ \colon\ x_{a^i_j}\ \longmapsto\ y^{(i)}_j\,,
\]
where $A_i\colon a^i_1<a^i_2<\dotsb<a^i_{|A_i|}$ and 
${\bf Y}^{(i)}=(y^{(i)}_1,y^{(i)}_2,\dotsc)$.

Let $\Bdot$ be a set composition of $[n]$ and set
$\beta:=|\Bdot|$. 
Let
$\calS_\beta:=\calS_{\beta_1}\times\calS_{\beta_2}\times\dotsb\times\calS_{\beta_s}$,
where $s$ is the length of $\Bdot$.
Identify $\calS_\beta$ as a subgroup of $\calS_n$ where an element
$u=(u^1,u^2,\dotsc,u^s)\in\calS_\beta$ is identified with 
$u^1\times u^2\times\dotsb\times u^s\in\calS_n$.
The longest element in $\calS_\beta$ is
$\omega_\beta:=\omega_{\beta_1}\times\dotsb\times\omega_{\beta_s}$.
Define the permutation $\varepsilon(\Bdot)\in\calS_n$ by 
 \begin{eqnarray}\label{E:epsilon1}
  \varepsilon(\Bdot)|_{B_i}&=& e\ \in\ \calS_{\beta_i},\ \\
  \varepsilon(\Bdot) &\colon& B_i\ \to\ \{\beta_1+\dotsb+\beta_{i-1}+1,\,
              \dotsc,\, \beta_1+\dotsb+\beta_i\}\,,\label{E:epsilon2}
 \end{eqnarray}
where $\beta_0=0$ and $\beta_{s+1}=n$.
Set $\zeta(\Bdot):=\omega_n\cdot\omega_\beta\cdot\varepsilon(\Bdot)$
and $u^R\cdot\zeta(\Bdot)=u^s\times\dotsb\times u^2\times
u^1\cdot\zeta(\Bdot)$. 
When $s=2$, we are in the situation of Definition~\ref{D:varepsilon}, so that 
$\Bdot=(P,Q)$ with $\varepsilon(\Bdot)=\varepsilon(P,Q)$
and $\zeta(\Bdot)=\zeta(P,Q)$.

\begin{thm}\label{T:decomposition}
 Suppose that $\Adot$ is a set composition of\/ $[m]$ with $s$ parts.
 For $w\in\calS_m$ suppose that $\beta=(\beta_1,\dotsc,\beta_s)$ is a 
 composition of $n\geq m$ such that $\beta_i\geq|A_i|$
 for each $i=1,\dotsc,s$, and
\[
  \psi_{\Adot}(\calG_w)\ \in\ 
   R_{\beta_1}({\bf Y}^{(1)})\otimes R_{\beta_2}({\bf Y}^{(2)})
   \otimes\dotsb\otimes R_{\beta_s}({\bf Y}^{(s)})\,.
\]
 Then, for any set composition $\Bdot$ of $n$ with $|\Bdot|=\beta$ and
 $B_i\cap[m]=A_i$ for each $i=1,\dotsc,s$, we have the formula
\[
  \psi_{\Adot}(\calG_w)\ =\ 
   \sum_{u\in\calS_\beta} 
    c^{u\cdot\varepsilon(\Bdot)}_{w,\,\varepsilon(\Bdot)}\
    \calG_{u^1}({\bf Y}^{(1)})\cdot
    \calG_{u^2}({\bf Y}^{(2)})\cdots\calG_{u^s}({\bf Y}^{(s)})\,.
\]
 Similarly, if $\beta$ satisfies
\[ 
  \psi_{\Adot}(R_n(x))\ \subset\ 
   R_{\beta_1}({\bf Y}^{(1)})\otimes R_{\beta_2}({\bf Y}^{(2)})
   \otimes\dotsb\otimes R_{\beta_s}({\bf Y}^{(s)})\,,
\]
 and $\Bdot$ is as before, then we have the formula 
\[
  \psi_{\Adot}(\calH_w)\ =\ 
   \sum_{u\in\calS_\beta} 
    c^{\omega_n w}_{\varepsilon(\Bdot),\, 
          \omega_n\cdot u\cdot\varepsilon(\Bdot)}\
    \calH_{u^1}({\bf Y}^{(1)})\cdot
    \calH_{u^2}({\bf Y}^{(2)})\cdots\calH_{u^s}({\bf Y}^{(s)})\,.
\]
\end{thm}

We prove this using a general form of the pattern map of Section~\ref{S:four},
where the torus $T_0\simeq\mathbb{C}^\times$ acts on $\mathbb{C}^n$ with weight
spaces of dimensions $\alpha_1,\alpha_2,\dotsc,\alpha_s$.
Then the flag variety of the centralizer of this torus is the product
of flag manifolds
\[
  \calF^\alpha\ :=\ 
   \Fl \mathbb{C}^{\alpha_1}\ \times\ \Fl\mathbb{C}^{\alpha_2}\ \times\ 
    \dotsb\ \times\ \Fl\mathbb{C}^{\alpha_s}\,.
\]
Given a set composition $\Adot$ of $[m]$ with $|\Adot|=\alpha$, there is a section
\[
   \varphi_\Adot\ \colon\ \calF^\alpha\ \longrightarrow\ \calF^{T_0}\ 
    \subset\ \calF
\]
of the pattern map analogous to the section $\varphi_{P,Q}$
of Section~\ref{S:four}.
Its action on the Grothendieck rings of $\calF^\alpha$ and $\calF$ is also
analogous.
We state the analog of Proposition~\ref{P:pattern_map}.

\begin{prop}
  $\varphi^*_\Adot=\overline{\psi_\Adot}$, the map of Grothendieck rings induced
  by $\psi_\Adot$.
\end{prop}

We now determine the map $\varphi^*_\Adot$ in terms of the basis of Schubert
structure sheaves. 
This begins with a geometric fact that generalizes Lemma~4.5.1
of~\cite{BS98}.

\begin{lemma}\label{L:multiProd}
  Let $\Adot=(A_1,\dotsc,A_s)$ be a set partition of\/ $[m]$ and set 
  $\alpha:=|\Adot|$. 
  Then there are flags $\Gdot$ and $\Gpdot$ in general position such that 
  for any $u\in\calS_\alpha$, we have 
\begin{eqnarray*}
    \varphi_{\Adot}\bigl(X_{u^1}\times\dotsb\times X_{u^s}\bigr)
   &\subset& X_{u^R\cdot\zeta(\Adot)}\Gdot\\
  \varphi_{\Adot}(\calF^\alpha) &\subset& X_{\varepsilon(\Adot)}\Gpdot
\end{eqnarray*}
\end{lemma}

\begin{proof}
 We prove this by induction on the number of parts of $\Adot$.
 The main engine is Lemma~4.5.1 of~\cite{BS98}, which is also the initial case
 when $s=2$.
 Let $\Adot=(A_1,\dotsc,A_s)$ be a set composition of $[m+\alpha_s]$,
 where $\alpha_s=|A_s|$.
 Consider the composition $(P,Q)$ of $[m+\alpha_s]$ where $Q=A_s$ and
 $P=[m+\alpha_s]-A_s=A_1\cup\dotsb\cup A_{s-1}$. 
 Then the order preserving bijection $f\colon P\to [m]$ induces a set 
 composition $\Bdot:=f(\Adot)$ of $[m]$ with $s{-}1$ parts.
 Set $v:=u^1\times\dotsb\times u^{s-1}
          \in\calS_{\alpha_1}\times\dotsb\times\calS_{\alpha_{s-1}}$.
 By the induction hypothesis and the definitions of $\varphi_{P,Q}$,
 $\varphi_\Adot$, and $\varphi_\Bdot$, we have
 \begin{eqnarray*}
    \varphi_{\Adot}(X_{u^1}\times\dotsb\times X_{u^s})&=&
    \varphi_{P,Q}\bigl(\varphi_\Bdot(X_{u^1}\times\dotsb\times X_{u^{s-1}})
          \times X_{u^s}\bigr)\\
    &\subset&
 \varphi_{P,Q}\bigl(X_{v^R\cdot\zeta(\Bdot)}\times X_{u^s}\bigr)\\
      &\subset&
     X_{u^s\times\left(v^R\cdot\zeta(\Bdot)\right)\cdot \zeta(P,Q)}\,.
 \end{eqnarray*}
 An easy calculation shows that 
\[
   u^s\times\left(v^R\cdot\zeta(\Bdot)\right)\cdot \zeta(P,Q)
   \ =\ 
   u^R\cdot \zeta(\Adot)\,,
\]
 which completes the proof of the first inclusion.
 The second inclusion is similar.
\end{proof}

Lemma~\ref{L:multiProd} implies that 
\[
  \varphi_{\Adot}(X_{u^1}\times\dotsb\times X_{u^s})\ \subset\ 
    X_{u^R\cdot\zeta(\Adot)}\Gdot\ \cap\ X_{\varepsilon(\Adot)}\Gpdot\,.
\]
 A dimension calculation (involving the lengths of the permutations
 $\varepsilon(\Adot)$ and $u^R\cdot\zeta(\Adot)$) shows that both sides have
 the same dimension. 
 As they are irreducible, we obtain the analog of Proposition~\ref{P:BS}
 and ultimately of Corollary~\ref{C:pushForward} and Theorem~\ref{T:pushForward}.

\begin{cor}\label{C:push}
 Let $\Adot=(A_1,A_2,\dotsc,A_s)$ be a set composition of\/ $[m]$ 
 and set $\alpha:=|\Adot|$.
 Then for any $u\in\calS_\alpha$ we have 
\begin{eqnarray*}
  (\varphi_\Adot)_*\bigl( 
    [\calO_{X_{u^1}}]\otimes\dotsb\otimes[\calO_{X_{u^s}}]\bigr) &=&
    [\calO_{X_{u^R\cdot\zeta(\Adot)}}]\cdot[\calO_{X_{\varepsilon(\Adot)}}]\\
 (\varphi_\Adot)_*\bigl( 
    [\calI_{u^1}]\otimes\dotsb\otimes[\calI_{u^s}]\bigr)&=&
    [\calI_{u^R\cdot\zeta(\Adot)}]\cdot[\calI_{\varepsilon(\Adot)}]\,.
\end{eqnarray*}
\end{cor}

Theorem~\ref{T:decomposition} follows from Corollary~\ref{C:push} in the
same manner as Theorem~\ref{T:SubsComp} followed from the results for $s=2$. 
As noted in Remark~\ref{R:identities}, the different possible choices of $\beta$
and $\Bdot$ of Theorem~\ref{T:decomposition} lead to a plethora of identities
among certain Schubert structure constants.
We concentrate instead on a special case which is used in~\cite{BSY} to show
that the quiver constants of Buch and Fulton are Schubert structure
constants.\smallskip  

Suppose first that the set composition $\Adot$ of $[m]$  
is sorted in that each $A_i$ is the interval
\[
  A_i\ =\ (a_{i-1},a_i]\ :=\ \{a_{i-1}+1,\, a_{i-1}+2,\,\dotsc,\, a_i\}\,,
\]
where $0=a_0<a_1<\dotsb<a_s=m$.
Then the specialization $\psi_\Adot(\calG_w)$ writes $\calG_w$ in terms of
Grothendieck polynomials in intervals of its variables.
Since the exponent of the variable $x_i$ in a monomial in $\calG_w$ for
$w\in\calS_m$ is at most $m{-}i$, we have
\[
  \psi_\Adot(\calG_w)\ \in\ 
   R_m({\bf Y}^{(1)})\otimes R_{m-a_1}({\bf Y}^{(2)})
   \otimes\dotsb\otimes  R_{m-a_{s-1}}({\bf Y}^{(s)})\,.
\]

We examine one rather simple choice for $\Bdot$ allowed by
Theorem~\ref{T:decomposition}.
Set $b_i:=(i{+}1)m-a_1-\dotsb-a_i$.
Define the set composition $\Apdot$ of $b_s$ by
\[
  A'_i\ :=\ (a_{i-1},a_i] \cup (b_{i-1}, b_i]\,.
\]
Then $A'_i\cap[m]=A_i$.
Also, since $a_s=m$, we have $b_s=b_{s-1}$, and so $A'_s=A_s$.
This choice of $\Apdot$ satisfies the hypotheses of
Theorem~\ref{T:decomposition}.

\begin{ex}\label{Ex:KMS}
 For example, suppose that $\Adot=
 (\{1,2\},\{3,4,5\},\{6,7,8,9\},\{10,11,12\})$.
 Then $a=(0,2,5,9,12)$, 
 $(b_0,b_1,b_2,b_3)=(12,22,29,32)$, and $\Apdot$ is the set
 composition of $\{1,\dotsc,32\}$ such that 
 \begin{eqnarray*}
  A'_1&=& \{1,2,\ 13, 14, \dotsc, 22\}\\
  A'_2&=& \{3, 4, 5,\  23, 24, \dotsc, 29\}\\
  A'_3&=& \{6, 7, 8, 9,\ 30, 31,32\}\\
  A'_4&=& \{10, 11, 12\} 
 \end{eqnarray*}
\end{ex}

A permutation $w\in\calS_m$ is {\bf Grassmannian} with descent at $k$ if
$w_i<w_{i+1}$, for all $i$, except possibly when $i=k$.
Associated to such a Grassmannian permutation $w\in\calS_m$ is a weakly
decreasing sequence, the partition 
\[
   (w_k-k, w_{k-1}-k+1, \dotsc, w_2-2, w_1-1)\,,
\] 
called its {\bf shape}.
If $w$ is Grassmannian with shape $\lambda$ and descent $k$, then we write
$w=w(\lambda,k)$. 
In the statements below, we use an exponent to indicate repeated parts in a
partition.
For example, if $\Apdot$ is the composition of Example~\ref{Ex:KMS}, then 
$\varepsilon(\Apdot)$ is Grassmannian with descent at $12$:
\[
  \varepsilon(\Apdot)\ =\ 
   (1, 2, \, 13, 14, 15, \, 23, 24, 25, 26,\, 30, 31, 32, \
    3, 4, \dotsc,12,\,16,\dotsc,22,\,27,28,29)\,,
\]
and it has shape $(20,20,20,17,17,17,17,10,10,10,0,0)=(20^3, 17^4, 10^3, 0^2)$.
In general, the permutation $\varepsilon(\Apdot)$ is Grassmannian.

\begin{lemma}\label{L:epsilonAdot}
  Let $\Apdot$ be as defined above.
  Then the permutation $\varepsilon(\Apdot)$ is Grassmannian with descent
  at $m$ and shape 
   $(\gamma_1^{\alpha_s},\, \gamma_2^{\alpha_{s-1}},\, \dotsc,\,
   \gamma_{s-1}^{\alpha_2},\, 0^{\alpha_1})$, 
 where $\gamma_i=b_{s-i}-a_{s-i}$ and $\alpha_i=a_i-a_{i-1}$.
\end{lemma}

We state the result of Theorem~\ref{T:decomposition} for specializations of this form.

\begin{thm}\label{T:Quiver}
  Let $0=a_0<a_1<\dotsb<a_s=m$ be integers and suppose that $\Adot$ is the set
  composition of $[m]$ such that $A_i=(a_{i-1}, a_i]$, for $i=1,\dotsc,s$.
  Let $b_i:=(i+1)m-a_1-\dotsb-a_i$ and set 
  $A'_i:=(a_{i-1}, a_i]\cup(b_{i-1}, b_i]$, for each $i=1,\dotsc,s$.
  Set $\beta=|\Apdot|$, so that $\beta_i=b_i+a_i-b_{i-1}-a_{i-1}$.
  Then, for any $w\in\calS_m$, we have
 \begin{equation}\label{E:Quiver}
  \psi_\Adot(\calG_w)\ =\ \sum_{u\in\calS_\beta}
     c^{u\cdot\varepsilon(\Apdot)}_{w,\,\varepsilon(\Apdot)}\ 
     \calG_{u^1}({\bf Y}^{(1)})\cdot
    \calG_{u^2}({\bf Y}^{(2)})\cdots\calG_{u^s}({\bf Y}^{(s)})\,.
 \end{equation}
\end{thm}

\begin{rem}
 The result of Theorem~\ref{T:Quiver}, like all such decomposition formulas in
 this paper for Grothendieck polynomials, gives an identical formula for
 Schubert polynomials, where the structure constants are those for the
 multiplication in the Schubert basis of cohomology.
 (These are the constants $c^w_{u,v}$ in this paper under the restriction
 $\ell(u)+\ell(v)=\ell(w)$.) 
 While many of these formulas were given in~\cite{BS98} in a more cryptic form,
 the formulas derived from Theorems~\ref{T:decomposition} and~\ref{T:Quiver} are  
 new.
 In the case of Theorem~\ref{T:Quiver}, the insertion algorithm of Bergeron and
 Billey~\cite{BeBi93} gives a formula for the constants 
 $c^{u\cdot \varepsilon(\Apdot)}_{w,\,\varepsilon(\Apdot)}$
 when 
\[
  \ell(w)+\ell(\varepsilon(\Apdot))\ =\
  \ell(u\cdot\varepsilon(\Apdot))
  \ =\ \ell(u) + \ell(\varepsilon(\Apdot))\,,
\]
 this last inequality is by an extension of Lemma~\ref{L:interval} (c).

 More specifically, Bergeron and  Billey devised an insertion algorithm 
 for $rc$-graphs to give a combinatorial proof of Monk's formula for Schubert
 polynomials. 
 While this insertion algorithm does not work in general to give a
 multiplication formula for Schubert polynomials, Kogan~\cite{Ko01,Ko03} 
 shows that it  does work to compute Schubert structure constants
 $c^w_{u,\,v}$ when the   permutations $u,v\in\calS_n$ have a
 special form:  $u$ has no descents after a position $m$ and $v$ is
 Grassmannian with descent at $m$.

 When $\ell(w)=\ell(u)$, this is exactly the situation for the constants  
 $c^{u\cdot \varepsilon(\Apdot)}_{w,\,\varepsilon(\Apdot)}$ of
 Theorem~\ref{T:Quiver}.
 Indeed, note that $w\in\calS_m$, so that $w$ has no descents after position
 $m$, and by Lemma~\ref{L:epsilonAdot}, $\varepsilon(\Apdot)$ is Grassmannian
 with descent at $m$.
 We deduce the following.

\begin{cor}
  Let $m$, $\Adot$, $\Apdot$, and $\beta$ be as in Theorem~$\ref{T:Quiver}$.
  Then the insertion algorithm of Bergeron and Billey computes the constants 
  $c^{u\cdot \varepsilon(\Apdot)}_{w,\,\varepsilon(\Apdot)}$
   of~\eqref{E:Quiver}, when $\ell(u)=\ell(w)$.  
\end{cor}

\end{rem}

A very special case of Theorem~\ref{T:Quiver} is when the sequence 
$a=(1,2,\dotsc,m)$.
Then, for $w\in\calS_m$, $\psi_{\Adot}(\calG_w)$ writes
the Grothendieck polynomial as a sum of monomials.
In this case, $\Adot$ consists of singletons with $A_i=\{i\}$
and $b_i=(i+1)m-\binom{i+1}{2}$, so that $\Apdot$ is a set composition of
$\binom{m+1}{2}$ with  
\[
   A'_i\ =\ \{i,\ im-{\textstyle\binom{i}{2}}+1,\,
                  im-{\textstyle\binom{i}{2}}+2,\,\dotsc,\,
                  im-{\textstyle\binom{i}{2}}=
              (i+1)m-{\textstyle\binom{i+1}{2}} \}\,,
\]
and $\varepsilon(\Apdot)$ is the Grassmannian permutation
with descent at $m$ and whose shape is
 \begin{equation}\label{E:lambda_def}
   \lambda\ =\ ({\textstyle\binom{m}{2}},\,\dotsc,\,
      (i{-}1)m-{\textstyle\binom{i}{2}},\,\dotsc,\,2m{-}3,\,m{-}1,\,0)\,.
 \end{equation}
Since $\varepsilon(\Apdot)$ is Grassmannian of descent $m$, 
$c^{u\cdot \varepsilon(\Apdot)}_{w,\,\varepsilon(\Apdot)}=0$ unless 
$u\cdot \varepsilon(\Apdot)\geq_m w$.
As $w$ has no descents after $m$, the characterization of the $m$-Bruhat order
in~\cite{BS98} implies that $u\cdot \varepsilon(\Apdot)$ also has no descents
after position $m$.
The Grothendieck polynomials $\calG_{u^i}({\bf Y}^{(i)})$ are pure
powers of the first variable in ${\bf Y}^{(i)}$, which forces
$u^i\in\calS_{m+1-i}$ to have the form 
$[l_i,1,\dotsc,\widehat{l_i},\dotsc,m+1-i]$.
Lastly, in $u\cdot\varepsilon(\Apdot)$, the component $u^i$ of $u$ permutes the
values in the interval 
$\Big((i-1)m-\binom{i-1}{2},\, im-\binom{i}{2}\Big]=A'_i\cap [m]$.
Thus
\[
   u\cdot\varepsilon(\Apdot)\ =\ 
   [l_1,\, l_2+m,\, l_3+2m-1,\, l_4+3m-3,\,
    \dotsc,\,l_m+(m{-}1)m-{\textstyle\binom{m-1}{2}},\,\dotsc\ ]\,,
\]
the remaining values taken in increasing order.
But this permutation is Grassmannian with descent at $m$ and shape
\begin{equation}\label{E:lambdaPrime}
     \lambda+(l_m,\dotsc,l_2,l_1)\,,
\end{equation}
which is still a partition as 
$\lambda_{m-1}-\lambda_{m+1-i}=m-i\geq l_i$.

Since $\frakS_u(x)=x_1^{l-1}$ for $u=[l,1,2,\dotsc,\widehat{l},\dotsc]$, we
deduce the following corollary.

\begin{cor}
  Let $w\in \calS_n$.
  The coefficient of $x_1^{l_1} x_2^{l_2}\dotsb x_m^{l_m}$ in the Schubert
  polynomial $\frakS_w(x)$ is the Schubert structure constant
  $c^{w(\lambda',m)}_{w, w(\lambda)}$, where $\lambda$ is defined
  by~\eqref{E:lambda_def} and $\lambda'$ is defined by~\eqref{E:lambdaPrime}. 
\end{cor}

 In~\cite{BS02} the coefficient of a monomial appearing in a Schubert polynomial
 was shown to be the intersection number of that Schubert class with a product
 of special Schubert classes corresponding to that monomial.
 This Corollary shows that the coefficient of a monomial in a
 Schubert polynomial is given by intersecting with a Schubert
 variety corresponding to that monomial.

\bibliographystyle{amsplain}

\providecommand{\bysame}{\leavevmode\hbox to3em{\hrulefill}\thinspace}
\providecommand{\MR}{\relax\ifhmode\unskip\space\fi MR }
\providecommand{\MRhref}[2]{%
  \href{http://www.ams.org/mathscinet-getitem?mr=#1}{#2}
}
\providecommand{\href}[2]{#2}

\end{document}